% !TEX spellcheck = English
\magnification=1000
\def\arXiv{oui}
  % MISE EN PAGE 
  
 \def\oui{oui} 
  
  %\magnification 1100
\ifx\arXiv\oui
\else
 \pdfpagewidth=210truemm
 \pdfpageheight=297truemm 
\fi
  
  %
  % Definitions pour les notes en bas de page
  % [Format : \note{texte de la note} avec numerotation automatique]
  %

  %
  %\anote=\footnote de plain (2 arguments)
  %

  \catcode`@=12 % at signs are no longer letters

  %reference de note
 \def\defrefnote#1{\definexref{#1}{{\the\footnotenumber}}{refnotes}}

  %
  % Fin notes
  %

  %\def\eqlgno#1{$$\leqalignno{#1}$$}

\ifx\couleurs\oui
\input graphicx
 \pdfpagewidth=210truemm
 \pdfpageheight=297truemm 
 \voffset=-5mm
\fi

\input eplain.tex
\expandafter\def\expandafter\newdimen\expandafter{\newdimen}

% couleurs
% Attention : \illust et graphicx nŽcessitent \couleurs={oui}
\ifx\couleurs\oui
%\expandafter\def\expandafter\newdimen\expandafter{\newdimen}
\beginpackages
\usepackage{color}
\endpackages 
 \pdfpagewidth=210truemm
 \pdfpageheight=297truemm 
\long\def\rge#1{{\color{red}#1}}

\definecolor{bleu-iecn}{cmyk}{.98,.13,.1,.55}

\else
\long\def\rge#1{#1}

\fi

\makeatletter
\def\numberedfootnote{%
ÊÊ\global\advance\footnotenumber by 1
ÊÊ\@eplainfootnote{{\number\footnotenumber}}%
}%
\def\makecolumns#1/#2 {\par \begingroup
ÊÊ \@columndepth = #1
ÊÊ \advance\@columndepth by -1
ÊÊ \divide \@columndepth by #2
ÊÊ \advance\@columndepth by 1
ÊÊ \@linestogoincolumn = \@columndepth
ÊÊ \@linestogo = #1
ÊÊ \currentcolumn = 1
ÊÊ \def\@endcolumnactions{%
ÊÊÊÊÊÊ\ifnum \@linestogo<2
ÊÊÊÊÊÊÊÊ \the\crtok \egroup \endgroup \par % End \valign and \makecolumns.
ÊÊÊÊÊÊ\else
ÊÊÊÊÊÊÊÊ \global\advance\@linestogo by -1
ÊÊÊÊÊÊÊÊ \ifnum\@linestogoincolumn<2
ÊÊÊÊÊÊÊÊÊÊÊÊ\global\advance\currentcolumn by 1
ÊÊÊÊÊÊÊÊÊÊÊÊ\global\@linestogoincolumn = \@columndepth
ÊÊÊÊÊÊÊÊÊÊÊÊ\the\crtok
ÊÊÊÊÊÊÊÊ \else
ÊÊÊÊÊÊÊÊÊÊÊÊ&\global\advance\@linestogoincolumn by -1
ÊÊÊÊÊÊÊÊ \fi
ÊÊÊÊÊÊ\fi
ÊÊ }%
ÊÊ \makeactive\^^M
ÊÊ \letreturn \@endcolumnactions
ÊÊ \@columnwidth = \hsize
ÊÊÊÊ \advance\@columnwidth by -\parindent
ÊÊÊÊ \divide\@columnwidth by #2
ÊÊ \penalty\abovecolumnspenalty
ÊÊ \noindent % It's not a paragraph (usually).
ÊÊ \valign\bgroup
ÊÊÊÊ &\hbox to \@columnwidth{\strut \hsize = \@columnwidth ##\hfil}\cr
}%
\makeatother

\lefteqnumbers
% choix numerotation gauche ou droite: defaut=gauche, \numadroite{oui} -> a droite
   \def\testd{oui}
   \def\choixlat{\ifx\numadroite\testd\righteqnumbers
            \else  \lefteqnumbers\fi}
    \choixlat

% modif format note
\catcode`@=\letter
\def\@eplainfootnote#1{\let\@sf\empty % parameter #2 (the text) is read later
  \ifhmode\edef\@sf{\spacefactor\the\spacefactor}\/\fi
  \global\advance\hlfootlabelnumber by 1
  \hlstart@impl{foot}{\hlfootlabel}%
  \hldest@impl{footback}{\hlfootbacklabel}%
  \hbox{$^{(#1)}$}%
  \hlend@impl{foot}%
  \@sf\vfootnote{#1.}%
}%
\catcode`@=\other

  \interfootnoteskip=0pt
  \let\note=\numberedfootnote
  \everyfootnote={\eightpoint\leftskip=5truemm\rightskip5truemm}
  
  \hsize150truemm\vsize 240truemm\hoffset=5truemm

  \pretolerance=500\tolerance=1000\brokenpenalty=5000
  \parindent3mm
  
  \countdef\temps=170
  \temps=\time
  \countdef\nminutes=171{\nminutes = \time}
  \countdef\nheure=172
  \def\heure{\begingroup                     % heure a la francaise
     \temps = \time \divide\temps by 60
     \nheure = \temps                        % l'heure, de 0 \`a 23
     \nminutes = \time
     \multiply\temps by 60
     \advance\nminutes by -\temps            % Les minutes, de 0 a 59.
     \ifnum\nminutes<10 \toks1 = {0}%
     \else\toks1 = {}%
     \fi
     \number\nheure h\the\toks1 \number\nminutes  
  \endgroup}%

  \newcount\chstart
  \chstart=\pageno
 \headline={\ifnum\pageno=\chstart {\hfill} \else {\hss \tenrm --\ \folio\ --\hss}\fi}
  \footline={\hfill}
  \normalbaselines
  \frenchspacing
    \def\dater{\vglue-10mm\rightline{(\the\day/\the\month/\the\year)}}
  \def\dateheure{\vglue-10mm\rightline{(\the\day/\the\month/\the\year,\ \heure)}}
 % justification ˆ droite

  \newif\ifpagetitre \pagetitretrue
\newtoks\hautpagetitre \hautpagetitre={\hfill}
\newtoks\baspagetitre \baspagetitre={\hfill}
\newtoks\auteurcourant \auteurcourant={\hfill}
\newtoks\titrecourant \titrecourant={\hfill}
\newtoks\hautpagegauche
\newtoks\hautpagedroite
\newtoks\hautpagemilieu
\hautpagemilieu={\tenrm\hfil -- \folio\ -- \hfil}
\hautpagegauche={\ifx\midfolio\oui\the\hautpagemilieu\else\tenrm\folio\hfill\the\auteurcourant\hfill\fi}
\hautpagedroite={\ifx\midfolio\oui\the\hautpagemilieu\else\hfill\the\titrecourant\hfill\tenrm\folio\fi}
\newtoks\baspagegauche \baspagegauche={\hfil}
\newtoks\baspagedroite \baspagedroite={\hfil}
\headline={\ifpagetitre\the\hautpagetitre
\else\ifodd\pageno\the\hautpagedroite\else\the\hautpagegauche\fi\fi }
\footline={\ifpagetitre\the\baspagetitre
\else\ifodd\pageno\the\baspagedroite
\else\the\baspagegauche\fi\fi \global\pagetitrefalse}

\def\pageblanche{\vfill\eject\pagetitretrue
\null\vfill\eject
\pagetitretrue
}
\def\chgtpage{\ifodd\pageno \else
\pageblanche \fi \pagetitretrue\titreun=0\footnotenumber=0}

\def\chgtpageincrtitreun{\ifodd\pageno \else
\pageblanche \fi \pagetitretrue\footnotenumber=0}

\def\majnombres{\ifodd\pageno \else
\pageblanche \fi \pagetitretrue\hautpoly\titreun=0\footnotenumber=0}

\def\hautspages#1#2{\auteurcourant={\ninepcap#1}\titrecourant={\nineit#2}}

\ifnum\chstart=\pageno \pagetitretrue\fi
%\ifnum\pageno=\chstartÊ\pagetitretrue \fi
  
  %pour eviter les rectangles noirs:
  %\overfullrule=Opt
  % pour que les @ soient des lettres

  %% pour aller a la ligne dans un \proclaim:
  \def\PAR{\par}
  
  %% Note en marge gauche
  \def\leftnote#1{\vadjust{\setbox1=\vtop{\hsize 20mm\parindent=0pt\eightpoint
  \baselineskip=9pt\rightskip=4mm plus 4mm\vskip-4.7mm#1}\hbox{\kern-2cm\smash{\box1}}}}

  %% Pour couper automatiquement un titre trop long
  
  \def\raggedcenter{\leftskip=20pt plus 10em  % reglages initial 4em
       \rightskip=\leftskip 
        \parfillskip=0pt 
         \spaceskip=.3333em \xspaceskip=.5em 
          \pretolerance=9999 \tolerance=9999
           \hyphenpenalty=9999 \exhyphenpenalty=9999 }
           
  \def\titrecentre#1{{\parindent0mm\raggedcenter
       \spaceskip=.6em plus .2em minus .2em\xspaceskip=.6em plus .2em minus .2em
        \tit#1\par}}
        
% entourer un texte par un filet

  %%%%%% D\'ebut num\'erotation automatique des paragraphes et/ou des enonces

  \def\oui{oui}
  
\def\fontetitreun{\ifx\paradouze\oui\douzepts\gpdouze\twelvebf\textfont1=\twelveib\else
\quatorzepts\gpquatorze\fourteenbf\fi}

\def\fontetitreunl{\douzepts\textfont1=\twelveib\scriptfont1=\tenib\fourteenti}
 
 \def\fontetitredeux{\textfont1=\eleveni\ifx\paradouze\oui\onzepts\scriptfont1=\ninei\elevenit\else
                        \douzepts\twelveit\fi}
 
   \def\fontetitredeuxb{\ifx\paradouze\oui\onzepts\eleventi\gponze\textfont1=\elevenib\scriptfont1=\nineib
                         \else\douzepts\twelveti\scriptfont1=\twelveib\scriptfont1=\tenib\gpdouze\fi}
                         
\def\fontetitredeuxl{\onzepts\textfont1=\elevenbf\scriptfont1=\ninebf\twelvebf}
  
\def\fontetitretrois{\textfont0=\elevenrm\scriptfont0=\eightrm\textfont1=\eleveni
                      \scriptfont1=\eighti\scriptscriptfont1=\sixi\elevenit}
                      
\def\fontetitrequatre{\textfont0=\elevenrm\scriptfont0=\eightrm\textfont1=\eleveni
                      \scriptfont1=\eighti\scriptscriptfont1=\sixi\elevenrm}
  
  \newcount\titreun\titreun=0
  \newcount\titredeux\titredeux=0
  \newcount\titretrois\titretrois=0
  \newcount\titrequatre\titrequatre=0
  \newcount\enonce\enonce=0
  
  \def\incr#1{\global\advance#1 by 1 {\the #1}}
  \def\avance#1{\global\advance#1 by 1}
  \def\init#1{\global#1=0}
  
  \long\def\Indentation#1#2{\setbox10=\hbox{\fontetitreun#1}
  	                    \ifdim\wd10 < 4mm
                         \setbox10=\hbox to 4mm{\box10\hfill}
                       \else\ifdim\wd10 < 6mm
                         \setbox10=\hbox to 6mm{\box10\hfill}
  	                    \else\ifdim\wd10 < 8mm
                         \setbox10=\hbox to 8mm{\box10\hfill}
                       \else\ifdim\wd10 < 12mm
                         \setbox10=\hbox to 12mm{\box10\hfill}
                       \fi\fi\fi\fi
                       \dimen10=\hsize
                       \advance \dimen10 by -\wd10
                       \noindent \box10 %
                       \ignorespaces
                       \hbox{\vtop{\hsize=\dimen10\raggedright\noindent\fontetitreun#2}}}

  \long\def\paraun#1{\removelastskip\par\medskip\goodbreak\vskip0pt plus.01\vsize\penalty-100
                \vskip0pt plus-.01\vsize
  	              \init{\titredeux}\ifnum\optionparag=1{\init\eqnumber\init\enonce}\else{}\fi
                  \goodbreak{\fontetitreun
  	                \Indentation{\incr{\titreun}.\ }{\fontetitreun #1\par}}\nobreak\medskip}

 %
 % titres de paragraphes centres
 %
 \long\def\paraunc#1{\removelastskip\par\bigskip\goodbreak\vskip0pt plus.01\vsize\penalty-100
                \vskip0pt plus-.01\vsize
  	              \init{\titredeux}
                 \ifnum\optionparag=1{\init{\eqnumber}\init\enonce}\else{}\fi
                  \goodbreak
  	                {\parindent0mm\raggedcenter\fontetitreun\incr{\titreun}.\ 
                     \fontetitreun #1\par}\nobreak\medskip}
                     
% paragraphes de niveau 1 numŽrotŽs avec des chiffres romains
\newtoks\titreunl
\titreunl={\ifnum\titreun=1{I}\fi%
\ifnum\titreun=2{II}\fi%
\ifnum\titreun=3{III}\fi%
\ifnum\titreun=4{IV}\fi%
\ifnum\titreun=5{V}\fi%
\ifnum\titreun=6{VI}\fi%
\ifnum\titreun=7{VII}\fi%
\ifnum\titreun=8{VIII}\fi%
\ifnum\titreun=9{IX}\fi%
\ifnum\titreun=10{X}\fi%
\ifnum\titreun=11{XI}\fi%
\ifnum\titreun=12{XII}\fi%
\ifnum\titreun=13{XIII}\fi%
}
\long\def\paraunl#1{\removelastskip\par\bigskip\bigskip\goodbreak\vskip0pt plus.01\vsize\penalty-100
                \vskip0pt plus-.01\vsize
  	              \init{\titredeux}\ifnum\optionparag=1{\init\eqnumber\init\enonce}\else{}\fi
                  \goodbreak{\fontetitreunl
  	                \Indentation{\global\advance\titreun by 1{\the\titreunl}.\ }{\fontetitreunl #1\par}}\nobreak\smallskip}

% paragraphes de niveaux 2
  
  \long\def\paradeux#1{\init{\titretrois}\vskip0pt plus.01\vsize\penalty-10
                \vskip0pt plus-.01\vsize\ifx \elie\oui\medskip\ifnum\titredeux>0\medskip\fi\fi
                 \Indentation{\fontetitredeux\the\titreun${\cdot}$\incr{\titredeux}.}
                              {\fontetitredeux\textfont1=\eleveni#1}\nobreak\par }
  
  \long\def\paradeuxb#1{\init{\titretrois}\vskip0pt plus.001\vsize\penalty-10
                \vskip0pt plus-.01\vsize{\ifx \elie\oui\medskip\ifnum\titredeux>0\medskip\fi\fi
                  \Indentation
  {\fontetitredeuxb\the\titreun${\cdot}$\incr{\titredeux}.}{ \fontetitredeuxb#1}}\nobreak
\smallskip}

% paragraphes de niveau 2 numŽrotŽs avec des lettres capitales 
\newtoks\titredeuxl
\titredeuxl={\ifnum\titredeux=1{A}\fi%
\ifnum\titredeux=2{B}\fi%
\ifnum\titredeux=3{C}\fi%
\ifnum\titredeux=4{D}\fi%
\ifnum\titredeux=5{E}\fi%
\ifnum\titredeux=6{F}\fi%
\ifnum\titredeux=7{G}\fi%
\ifnum\titredeux=8{H}\fi%
\ifnum\titredeux=9{I}\fi%
\ifnum\titredeux=10{J}\fi%
\ifnum\titredeux=11{K}\fi%
\ifnum\titredeux=12{L}\fi%
\ifnum\titredeux=13{M}\fi%
}
 \long\def\paradeuxl#1{\init{\titretrois}\vskip0pt plus.001\vsize\penalty-10
                \vskip0pt plus-.01
                \vsize \bigskip%
                  \Indentation
     {\fontetitredeuxl\global\advance\titredeux by 1
  \quad \the\titreunl${\cdot}$\the\titredeuxl.}{ \fontetitredeuxl#1}
  \removelastskip\nobreak\smallskip}
  
% paragraphes de niveaux 3

  \long\def\paratrois#1{\init{\titrequatre}\ifdim\lastskip<\smallskipamount
                \removelastskip\smallskip\fi
                 \vskip0pt plus.01\vsize\penalty-10
                  \vskip0pt
plus-.01\vsize{\ifx \elie\oui\ifnum\titretrois>0\medskip\fi\fi
\Indentation{\fontetitretrois\the\titreun${\cdot}$\the\titredeux${\cdot}$\incr{\titretrois}.\ }
  {\hskip0mm\baselineskip=14pt\fontetitretrois#1}\nobreak\smallskip}}
  
 % paragraphes de niveau 3 numŽrotŽs avec des chiffres mais niveaux 1 et 2 chiffres romains et lettres 
  
  \long\def\paratroisl#1{\init{\titrequatre}\ifdim\lastskip<\smallskipamount
                \removelastskip\fi
                 \vskip0pt plus.01\vsize\penalty-10
                  \vskip0pt
plus-.01\vsize\ifx \elie\oui\bigskip
\fi
\Indentation{\fontetitretrois\quad \quad \the\titreunl{${\cdot}$}\the\titredeuxl${\cdot}$\incr{\titretrois}.\ }
  {\hskip0mm\fontetitretrois#1}\nobreak\smallskip}

% paragraphes de niveaux 4

  \long\def\paraquatre#1{\ifdim\lastskip<\smallskipamount
                \removelastskip\smallskip\fi
                 \vskip0pt plus.01\vsize\penalty-10
                  \vskip0pt
                  plus-.01\vsize\par
                %\bigskip
 
\Indentation{\fontetitrequatre \the\titreun{${\cdot}$}\the\titredeux${\cdot}$\the\titretrois${\cdot}$\incr{\titrequatre}.\ }
{\hskip0mm\fontetitrequatre#1}\nobreak\smallskip}

% paragraphes de niveau 4 numŽrotŽs avec des lettres minuscules

\newtoks\titrequatrel
\titrequatrel={\ifnum\titrequatre=1{a}\fi%
\ifnum\titrequatre=2{b}\fi%
\ifnum\titrequatre=3{c}\fi%
\ifnum\titrequatre=4{d}\fi%
\ifnum\titrequatre=5{e}\fi%
\ifnum\titrequatre=6{f}\fi%
\ifnum\titrequatre=7{g}\fi%
\ifnum\titrequatre=8{h}\fi%
\ifnum\titrequatre=9{i}\fi%
\ifnum\titrequatre=10{j}\fi%
\ifnum\titrequatre=11{k}\fi%
\ifnum\titrequatre=12{l}\fi%
\ifnum\titrequatre=13{m}\fi%
}
\long\def\paraquatrel#1{\ifdim\lastskip<\smallskipamount
                \removelastskip\smallskip\fi
                 \vskip0pt plus.01\vsize\penalty-10
                  \vskip0pt
                  plus-.01\vsize{\bigskip
\Indentation{\global\advance\titrequatre by 1
\fontetitrequatre\quad \quad \quad \the\titreunl${\cdot}$\the\titredeuxl${\cdot}$\the\titretrois${\cdot}$\the\titrequatrel.\ }
{\hskip0mm\fontetitrequatre#1}\nobreak\smallskip}}

% Pour memoriser le no de section d'un titre de paragraphe
\ifx\optionkeys\oui
\def\drefun#1{\definexref{¤#1}{{\the\titreun}}{}} 
\def\drefdeux#1{\definexref{¤#1}{{\the\titreun}.{\the\titredeux}}{}}
\def\dreftrois#1{\definexref{¤#1}{{\the\titreun}.{\the\titredeux}.{\the\titretrois}}{}}
\else
\def\drefun#1{\definexref{prg#1}{{\the\titreun}}{}} 
\def\drefdeux#1{\definexref{prg#1}{{\the\titreun}.{\the\titredeux}}{}}
\def\dreftrois#1{\definexref{prg#1}{{\the\titreun}.{\the\titredeux}.{\the\titretrois}}{}}
\fi

% Dans le cas de \drefun, coder \paraun{Titre}\drefun{label} dans le cas de \drefdeux, on peut coder
%                                \paradeux{Sous-titre\drefdeux{label}}
% pour appeler la reference "label", taper \refn{¤label}
%
% en cas de probl\`emes avec option_keys: coder \def\optionkeys{oui} et remplacer ¤ par prg dans \refn

%% num\'erotation des \'enonc\'es
%%Exemples
%%\propt{leb}{ (Lebesgue)}{ L'espace $L^1$ est complet.}
%%\par
%%Ca marche, comme le dit le \ref{leb}.
%%\Propt{leb2}{L'espace $L^1$ est vraiment complet.}\par
%%Ce \ref{leb2} est moins bien que les Th\'eor\`emes \ref{sleb} et \ref{sleb2}.

  \long\def\propdeux#1#2#3#4{%
       \avance{\enonce}
       \leavevmode\edef\temp{#2}%
         \ifx\temp\empty 
          \else
           \definexref{#2}{#1~{\the\titreun.\the\enonce}}{enonces}
            \definexref{s#2}{{\the\titreun.\the\enonce}}{enonces}
             \fi
\smallskip
      \noindent{\bf#1\ {\bf\the\titreun.\the\enonce{#3}.}\enspace}{\sl#4\par}%
      \ifdim\lastskip<\medskipamount \removelastskip\penalty55\par \fi
   }

  \long\def\propun#1#2#3#4{%
      \avance{\enonce}
       \leavevmode\edef\temp{#2}%
        \ifx\temp\empty 
          \else
           \definexref{#2}{#1~{\the\enonce}}{enonces}
            \definexref{{s#2}}{{\the\enonce}}{enonces}
             \fi
   \par 
     \noindent{\bf#1\ {\bf\the\enonce{#3}.}\enspace}{\sl#4\par}%
     \ifdim\lastskip<\medskipamount \removelastskip\penalty55\medskip\fi
  }
  
  \long\def\prop#1#2#3#4{\ifnum\optionparag=1
                          \propdeux{#1}{#2}{\textfont1=\elevenib#3}{#4} \else\propun{#1}{#2}{\textfont1=\elevenib#3}{#4}\fi}

  \long\def\propt#1#2#3{\ifx\tpf\oui \prop{Th\'eo\-r\`eme}{#1}{#2}{#3}\par
                       \else\prop{Theorem}{#1}{#2}{#3}\par\fi}
  \long\def\Propt#1#2{\propt{#1}{}{#2}}
  \long\def\propl#1#2#3{\ifx\tpf\oui\prop{Lem\-me}{#1}{#2}{#3}\par
                         \else\prop{Lemma}{#1}{#2}{#3}\par\fi}
  \long\def\Propl#1#2{\propl{#1}{}{#2}}
  \long\def\propc#1#2#3{\ifx\tpf\oui\prop{Corol\-laire}{#1}{#2}{#3}\par
                         \else\prop{Corollary}{#1}{#2}{#3}\par\fi}
  
  \long\def\propp#1#2#3{\prop{Pro\-po\-si\-tion}{#1}{#2}{#3}\par}
  \long\def\Propp#1#2{\propp{#1}{}{#2}} 
  \long\def\propd#1#2#3{\ifx\tpf\oui\prop{D\'efi\-nition}{#1}{#2}{#3}\par
                       \else\prop{Definition}{#1}{#2}{#3}\par\fi} 
  
  \long\def\proptd#1#2#3{\ifx\tpf\oui\prop{Th\'eor\`eme et d\'efi\-nition}{#1}{#2}{#3}\par
                       \else\prop{Theorem and definition}{#1}{#2}{#3}\par\fi}

  % Pour une numerotation des enonces et des formules avec un seul nombre 
  % coder \optionparag=2

  % Choisir numerotation manuelle ou automatique des paragraphes: 
  % utiliser \section au lieu de \paraun et choisir \optionparag=2
  % Cela induit aussi une num\'erotation des formules avec un seul nombre
  
  \newcount\optionparag\optionparag=1
  
  \long\def\section#1#2{\ifnum\optionparag=1 \paraun{#2} 
                        \else\goodbreak{\fontetitreun
  	                \Indentation{#1.\ }{#2}}\nobreak\smallskip\fi}

    % numerotation automatique (par defaut)
    % numerotation manuelle
  \def\eqconstruct#1{\ifnum\optionparag=1{\the\titreun\hbox{$\cdot$}#1}\else{#1}\fi}

  %%%%%%% Fin num\'erotation automatique des paragraphes
  
  %% num\'erotation bibliographie
  
  \def\numref{oui}  % numeroter la biblio par defaut avec bibtem
  
  \newcount\mesref\mesref=0 
  \def\defbib#1{\ifx\numref\oui\global\advance\mesref by 1 \definexref{#1}{{\the
                 \mesref}}{}\else\definexref{#1}{#1}{}\fi}
  \def\bibtem#1{\defbib{#1}\item{\citer{#1}}}
  \def\citer#1{[\ref{#1}]}
  \def\citeplus#1#2{[\ref{#1}; #2]}

  % FONTES & FAMILLES
  
  \font\seventeenmsa=msam10 at 17pt    % symboles d'AMSTEX
  \font\fourteenmsa=msam10 at 14pt
  \font\twelvemsa=msam10 at 12pt
  \font\tenmsa=msam10                 
  \font\ninemsa=msam10 at 9pt 
  \font\eightmsa=msam10 at 8pt 
  \font\sevenmsa=msam7 
  \font\sixmsa=msam10 at 6pt
  \font\fivemsa=msam5
  \newfam\msafam\textfont\msafam=\tenmsa\scriptfont\msafam=\sevenmsa\scriptscriptfont\msafam=\fivemsa
  
  \font\seventeenbb=msbm10 at 17pt     % Lettres evidees pour titres
  \font\fourteenbb=msbm10 at 14pt
  \font\twelvebb=msbm10 at 12pt
  \font\tenbb=msbm10                   %Lettres evidees
  \font\ninebb=msbm10 at 9pt 
  \font\eightbb=msbm10 at 8pt 
  \font\sevenbb=msbm7 
  \font\sixbb=msbm10 at 6pt
  \font\fivebb=msbm5 
  \newfam\bbfam\textfont\bbfam=\tenbb\scriptfont\bbfam=\sevenbb\scriptscriptfont\bbfam=\fivebb
  \def\bb{\fam\bbfam\tenbb}%

  \font\seventeenscaln=eusm10 at 17pt   % Lettres de ronde
  \font\twelvescaln=eusm10 at 12pt
  \font\tenscaln=eusm10                
  \font\ninescaln=eusm10 scaled 900
  \font\eightscaln=eusm10 scaled 800
  \font\sevenscaln=eusm10 scaled 700
  \font\sixscaln=eusm10 scaled 600
   
  \newfam\scalnfam\textfont\scalnfam=\tenscaln\scriptfont\scalnfam=\sevenscaln\scriptscriptfont\scalnfam=\sixscaln
  \def\scaln{\fam\scalnfam\tenscaln}%
  \def\scal{\scaln}
  
  \font\tenscalb=eusb10                % Lettres de ronde grasses

  \font\sevenscalb=eusb10 scaled 700

  \newfam\scalbfam\textfont\scalbfam=\tenscalb\scriptfont\scalbfam=\sevenscalb
  %
  
  %
  % Romain
  %
  \font\fourteenrm=cmr12 scaled 1200
  \font\elevenrm=cmr10 at 11pt
  \font\twelverm=cmr12
  \font\ninerm=cmr9
  \font\eightrm=cmr8      
  \font\sevenrm=cmr7
  \font\sixrm=cmr6

  \font\seventeenpcap=cmcsc10 at 17pt
  \font\tenpcap=cmcsc10                        % Petites capitales
  \font\ninepcap=cmcsc9
  \font\eightpcap=cmcsc8
  \font\sevenpcap=cmcsc10 scaled 700
  
  \newfam\pcapfam\textfont\pcapfam=\tenpcap\scriptfont\pcapfam=\sevenpcap
  \def\pcap{\fam\pcapfam\tenpcap}
  
               % Gras romain (boldface)
                % Titres
  \font\seventeenrm=cmbx12 scaled 1400

  \font\fourteenbf=cmbx10 scaled 1400
  
  \font\twelvebf=cmbx12
  \font\elevenbf=cmbx10 at 11pt
  \font\ninebf=cmbx9  
  \font\eightbf=cmbx8
  \font\sixbf=cmbx6
  
  \font\tengot=eufm10                           % Lettres gothiques
   
  \font\eightgot=eufm10 at 8truept 
  \font\sevengot=eufm7 
  \font\sixgot=eufm10 at 6 truept 
   
  \newfam\gotfam
  \textfont\gotfam=\tengot\scriptfont\gotfam=\sevengot\scriptscriptfont\gotfam=\sixgot
  \def\got{\fam\gotfam\tengot}%

  %% Pour les titres (168pt)
  
  \def\tit{%
  \textfont0=\seventeenrm\scriptfont0=\tenrm\def\rm{\fam0\seventeenrm}%
  \textfont1=\seventeenib\scriptfont1=\twelveib%
  \textfont2=\seventeensy\scriptfont2=\twelvesy\scriptscriptfont2=\ninesy
  \textfont3=\seventeenex
  \textfont\itfam=\seventeenti
  \def\it{\fam\itfam\seventeenti}%
  \textfont\bbfam=\seventeenbb \scriptfont\bbfam=\twelvebb
  \def\bb{\fam\bbfam\seventeenbb}%
  \textfont\msafam=\seventeenmsa\scriptfont\msafam=\twelvemsa
  \textfont\scalnfam=\seventeenscaln
  \def\pcap{\fam\pcapfam\seventeenpcap}
  \normalbaselineskip=25pt\normalbaselines\rm}

  % italiques grasses pour titres
  \font\seventeenti=cmbxti10 scaled 1680
  
  \font\fourteenti=cmbxti10 at 14pt
  
  \font\twelveti=cmbxti10 scaled 1200
  \font\eleventi=cmbxti10 at 11pt

  %
  % italiques
  %
  \font\twelveit=cmti12	
  \font\elevenit=cmti10 scaled 1100
  \font\nineit=cmti9
  \font\eightit=cmti8
  \font\sevenit=cmti7

  %
  % italique mathematique gras pour titres
  %
 
 \font\seventeenib=cmmib10 scaled 1680
  \font\fourteenib=cmmib10 scaled 1400
  \font\twelveib=cmmib10 scaled 1200
  \font\elevenib=cmmib10 scaled 1100
  \font\tenib=cmmib10
\font\eightib=cmmib10 scaled 800
  \font\nineib=cmmib10 scaled 900
\font\sevenib=cmmib10 scaled 700
\font\sixib=cmmib10 scaled 600
\font\fiveib=cmmib10 scaled 500

\ifx\ITAN\oui
\else
\innernewfam\cmmibfam
\textfont\cmmibfam=\tenib
\scriptfont\cmmibfam=\sevenib
\scriptscriptfont\cmmibfam=\fiveib
\def\ib{\fam\cmmibfam\tenib}
\fi

  %
  % italique mathematique
  %
  \font\twelvei=cmmi10 scaled 1200
  \font\eleveni=cmmi10 scaled 1100
  \font\ninei=cmmi9
  \font\eighti=cmmi8 
  \font\seveni=cmmi7 			                
  \font\sixi=cmmi6
  
  \font\ninesl=cmsl9                    % slanted 
  \font\eightsl=cmsl8 
  \font\sevensl=cmsl10 at 7pt

  \font\ninett=cmtt9                    % typewriter
  \font\eighttt=cmtt8
  \font\seventt=cmtt10 scaled 700

  \font\seventeensy=cmsy10 scaled 1680    % symboles pour titres
  \font\fourteensy=cmsy10 scaled 1400
  \font\twelvesy=cmsy10 scaled 1176
  
  \font\ninesy=cmsy9                      % symboles
  \font\eightsy=cmsy8
  \font\sixsy=cmsy6
  \font\seventeenex=cmex10 at 17pt
  \font\fourteenex=cmex10 at 14pt
  \font\twelveex=cmex10 at 12pt
  \font\nineex=cmex10 at 9pt
  \font\eightex=cmex10 at 8pt
  \font\sevenex=cmex10 at 7pt
  \font\sixex=cmex10 at 6pt
  \font\fiveex=cmex10 at 5pt
  
   %% Lettres grecques plombees              %% Lettres grecques plombees
   
  \font\fourteengp=cmmi10 at 14pt
  
  \font\twelvegp=cmmib10 at 12pt
  \font\elevengp=cmmib10 at 11pt
  \font\tengp=cmmib10                          
  \font\ninegp=cmmib10 at 9pt 
  \font\eightgp=cmmib8 
  \font\sevengp=cmmib7 
  \font\sixgp=cmmib6

%\innernewfam\gpfam
%\textfont1\gpfam=\tengp
%\scriptfont1\gpfam=\sevengp
%\scriptscriptfont1\gpfam=\fivegp
%\def\gp{\fam\gpfam\tengp}

\def\gp{\textfont0=\tenbf\scriptfont0=\sevenbf\scriptscriptfont0=\fivebf
\textfont1=\tengp\scriptfont1=\sevengp\scriptscriptfont1=\fivegp}
  \def\gponze{\textfont0=\elevenbf\scriptfont0=\eightbf\scriptscriptfont0=\sixbf
           \textfont1=\elevengp\scriptfont1=\eightgp\scriptscriptfont1=\sixgp}
  \def\gpdouze{\textfont0=\twelvebf\scriptfont0=\tenbf\scriptscriptfont0=\ninebf
           \textfont1=\twelvegp\scriptfont1=\tengp\scriptscriptfont1=\ninegp}        
  
 \def\gpquatorze{\textfont0=\fourteenbf\scriptfont0=\twelvebf\scriptscriptfont0=\elevenbf
           \textfont1=\fourteengp\scriptfont1=\twelvegp\scriptscriptfont1=\elevengp}

  % CARACTERES SPECIAUX
  
  \expandafter\chardef\csname pre amssym.def at\endcsname=\the\catcode`\@
  \catcode`\@=11
  \def\undefine#1{\let#1\undefined}
  \def\newsymbol#1#2#3#4#5{\let\next@\relax
   \ifnum#2=\@ne\let\next@\msafam@\else
   \ifnum#2=\tw@\let\next@\bbfam@\fi\fi
   \mathchardef#1="#3\next@#4#5}
  \def\mathhexbox@#1#2#3{\relax
   \ifmmode\mathpalette{}{\m@th\mathchar"#1#2#3}%
   \else\leavevmode\hbox{$\m@th\mathchar"#1#2#3$}\fi}
  \def\hexnumber@#1{\ifcase#1 0\or 1\or 2\or 3\or 4\or 5\or 6\or 7\or 8\or
   9\or A\or B\or C\or D\or E\or F\fi}
  
  \def\setboxz@h{\setbox\z@\hbox}
  \def\wdz@{\wd\z@}
  \def\boxz@{\box\z@}
  
  \edef\msafam@{\hexnumber@\msafam}
  \mathchardef\dabar@"0\msafam@39
  
  \edef\bbfam@{\hexnumber@\bbfam}
  \def\widehat#1{\setboxz@h{$\m@th#1$}%
   \ifdim\wdz@>\tw@ em\mathaccent"0\bbfam@5B{#1}%
   \else\mathaccent"0362{#1}\fi}
  \def\widetilde#1{\setboxz@h{$\m@th#1$}%
   \ifdim\wdz@>\tw@ em\mathaccent"0\bbfam@5D{#1}%
   \else\mathaccent"0365{#1}\fi}
  \newsymbol\leqq 1335          % superieur ou egal(=)
  \newsymbol\leqslant 1336
  \newsymbol\lessgtr 1337       % superieur ou inferieur
  \newsymbol\backprime 1038     % apostrophe de gauche a droite
  \newsymbol\risingdotseq 133A  % egal entre points (bas puis haut)
  \newsymbol\fallingdotseq 133B % egal entre points (haut puis bas)
  \newsymbol\succcurlyeq 133C   % superieur ou egal tordu
  \newsymbol\geqq 133D          % inferieur ou egal(=)
  \newsymbol\geqslant 133E
  \newsymbol\nmid 232D
  \newsymbol\nexists 2040
  \newsymbol\smallsetminus 2272
  \newsymbol\varnothing 203F 
  \catcode`\@=\active

  % typographie francaise ou anglaise
  \catcode`\@=11
  \newcount\typofr\typofr=1
  
  \catcode`\;=\active
  \def;{\ifnum\typofr=1\relax\ifhmode\ifdim\lastskip>\z@\unskip\fi
     \kern.2em\fi\string;\else\string;\fi}
  
  \catcode`\:=\active
  \def:{\ifnum\typofr=1\relax\ifhmode\ifdim\lastskip>\z@\unskip\fi
  \penalty\@M\ \fi\string:\else\string:\fi}
  
  \catcode`\!=\active
  \def!{\ifnum\typofr=1\relax\ifhmode\ifdim\lastskip>\z@\unskip\fi
     \kern.2em\fi\string!\else\string!\fi}
  
  \catcode`\?=\active
  \def?{\ifnum\typofr=1\relax\ifhmode\ifdim\lastskip>\z@\unskip\fi
     \kern.2em\fi\string?\else\string?\fi}

  \def\francais{\typofr=1\def\tpf{oui}}
  \def\anglais{\typofr=2\def\tpf{non}\def\english{oui}}
  \def\oui{oui}
  \francais
  
  \catcode`\@=12
  %\catcode`\@=\active
  
  %pour rayer du texte en rouge

\ifx\textures\oui
\def\raye #1|{\leavevmode\setbox1=\hbox{#1}%
\raise .5pt\hbox to \wd1{\xleaders\hbox{\rge{$ \sct / $}%
\kern 1pt}\hfill\kern -1pt }\kern -\wd1 \unhbox1\relax }

\def\barre #1|{\leavevmode\setbox1=\hbox{#1}%
\rlap{\Red\vrule height 2.4pt depth -1.2pt width \wd1}\Black \unhbox1\relax}
\else
\def\raye #1|{\leavevmode\setbox1=\hbox{#1}%
\raise .5pt\hbox to \wd1{\xleaders\hbox{\rge{$ \sct / $}%
\kern 1pt}\hfill\kern -1pt }\kern -\wd1 \unhbox1\relax }

\def\barre #1|{\leavevmode\setbox1=\hbox{#1}%
\rlap{\color{red}\vrule height 2.4pt depth -1.2pt width \wd1}\color{black} \unhbox1\relax}

\fi
  
  % FIN CARACTERES SPECIAUX

  % MACROS DIVERSES
  
  \def\og{\leavevmode\raise.24ex\hbox{$\scriptscriptstyle\langle\!\langle\>$}}    % guillemets ouvrants
  \def\fg{\leavevmode\raise.24ex\hbox{$\scriptscriptstyle\>\rangle\!\rangle$}}    % guillemets fermants

  \def\z{{\bb Z}}
  \def\r{{\bb R}}
  
  \def\N{{\bb N}}
  \def\Q{{\bb Q}}

  \def\K{{\scal K}}

  \def\M{{\scal M}}
  \def\n{{\scal N}}
  \def\O{{\scal O}}
  \def\P{{\scaln P}}
  
  \def\R{{\scal R}}
  \def\s{{\cal S}}

  \def\frac#1#2{{#1\over #2}}
  \def\di#1#2{\sct#1\atop{\sct#2}}
  \def\tri#1#2#3{{\sct#1\atop\sct#2}\atop\sct#3}

  \def\qedbox{$\rlap{$\sqcap$}\sqcup$}           % carr\'e blanc fin de d\'emo
  \def\qed{\nobreak\hfill\penalty250 \hbox{}\nobreak\hfill\qedbox\par }

  \def\¤{\S\thinspace}

  \def\¥{$\bullet$ }
  
  %SYMBOLES MATHS
  
  \def\e{{\rm e}}
  \def\mod{\mathop{\rm mod}\nolimits}
  \def\md#1#2{\equiv#1\,({\rm mod\,}#2)}
  \def\no#1{\Vert#1\Vert}

  \def\epsilon{\varepsilon}

  \def\phi{\varphi}
  \def\theta{\vartheta}
  \def\rho{\varrho}
  \def\dm{{\textstyle{1\over 2}}}
  \def\txt{\textstyle}
  
  \def\sct{\scriptstyle}
  \def\pf{\noi{\it Proof. }}
  \def\nid{\ifnum\typofr=1\par\noindent{\it D\'emonstration. }\else\pf\fi}
  \def\noi{\noindent}
  \def\rem{\ifnum\typofr=1\noi{\it Remarque.}\ \else\noi{\it Remark.}\ \fi}
  \def\rems{\ifnum\typofr=1\noi{\it Remarques.}\ \else\noi{\it Remarks.}\ \fi}

  \def\1{{\bf 1}}
  \def\|{\Vert}

  \def\le{\leqslant}\def\leq{\leqslant}
  \def\geq{\geqslant}
  \def\wh{\widehat}

  \newsymbol\subsetneqq 2324
  \newsymbol\subsetneq 2328
  %% operateurs
  
  \def\card{\mathop{\rm card}\nolimits}

  \def\log{\mathop{\rm log}\nolimits}
  \def\ft#1#2{{\txt{#1\over #2}}}

  %% Grande asterisque operateur

  %% Grande asterisque sur la ligne

\def\Vbs#1{\bigg|#1\bigg|}

  %fleches

  %% Lettres plombees 'pauvres'
  \def\pmb#1{\setbox0=\hbox{#1}%
  \kern-.025em\copy0\kern-\wd0\kern.05em\copy0\kern-\wd0\kern-.025em\raise .0433em\box0 }

  % NOTES EN BAS DE PAGE
  
  % Pour que les accents se placent correctement en mode math en corps 8 et 6
  \skewchar\eighti='177 \skewchar\sixi='177
  \skewchar\eightsy='60 \skewchar\sixsy='60
  
  \def\eightpoint{%
  \textfont0=\eightrm\scriptfont0=\sixrm\scriptscriptfont0=\fiverm
  \def\rm{\fam0\eightrm}%
  \textfont1=\eighti\scriptfont1=\sixi
  \scriptscriptfont1=\fivei\def\oldstyle{\fam1\seveni}%
  \textfont2=\eightsy\scriptfont2=\sixsy\scriptscriptfont2=\fivesy
  \textfont3=\eightex\scriptfont3=\sixex
  \textfont\itfam=\eightit
  \def\it{\fam\itfam\eightit}%
  \textfont\slfam=\eightsl
  \def\sl{\fam\slfam\eightsl}%
  \textfont\bbfam=\eightbb \scriptfont\bbfam=\sixbb\scriptscriptfont\bbfam=\fivebb
  \def\bb{\fam\bbfam\eightbb}%
  \textfont\msafam=\eightmsa\scriptfont\msafam=\sixmsa
  \textfont\scalnfam=\eightscaln
  \def\scaln{\fam\scalnfam\eightscaln}
  \textfont\ttfam=\eighttt
  \def\tt{\fam\ttfam\eighttt}%
\textfont\gotfam=\eightgot
  \textfont\bffam=\eightbf\scriptfont\bffam=\sixbf\scriptscriptfont\bffam=\fivebf
  \def\bf{\fam\bffam\eightbf}%
  \ifx\ITAN\oui\else\textfont\cmmibfam=\eightib
       \scriptfont\cmmibfam=\sixib
        \scriptscriptfont\cmmibfam=\fiveib
         \def\ib{\fam\cmmibfam\eightib}
   \fi
  \textfont\pcapfam=\eightpcap
  \def\pcap{\fam\pcapfam\eightpcap}
  \abovedisplayskip=2pt plus2pt minus 2pt
  \belowdisplayskip=2pt plus1pt minus 2pt
  \abovedisplayshortskip= 1pt plus 2pt minus 1pt
  \belowdisplayshortskip= 1pt plus 2pt minus 1pt
  \smallskipamount=2pt plus 1pt minus 2pt
  \medskipamount=3pt plus 2pt minus 2pt
  \bigskipamount=7pt plus 3pt minus 3pt
  \setbox\strutbox=\hbox{\vrule height 5pt depth 2pt width 0pt}%
  \normalbaselineskip=9pt\normalbaselines\rm}

  \def\({\left(}
  \def\){\right)}
  
  \def\footnoterule{\kern -2pt\hrule width 7truecm\kern 2.4pt}
  
  %reference croisee a des numeros de notes
  \def\xnotedef#1{\definexref{#1}{\noexpand\number\footnotenumber}{Note}}%

  %Fin des definitions pour les notes en bas de page
  
  % PARAGRAPHES EN NEUF POINTS
  
  \def\ninepoint{%
  \textfont0=\ninerm\scriptfont0=\sixrm\scriptscriptfont0=\fiverm
  \def\rm{\fam0\ninerm}%
  \textfont1=\ninei\scriptfont1=\sixi
  \scriptscriptfont1=\fivei\def\oldstyle{\fam1\ninei}%
  \textfont2=\ninesy\scriptfont2=\sixsy\scriptscriptfont2=\fivesy
  \textfont3=\nineex\scriptfont3=\sixex
  \textfont\itfam=\nineit
  \def\it{\fam\itfam\nineit}%
  \textfont\slfam=\ninesl
  \def\sl{\fam\slfam\ninesl}%
  \textfont\bbfam=\ninebb\scriptfont\bbfam=\sixbb\scriptscriptfont\bbfam=\fivebb
  \def\bb{\fam\bbfam\ninebb}%
  \textfont\msafam=\ninemsa\scriptfont\msafam=\sixmsa\scriptscriptfont\msafam=\fivemsa
  \textfont\scalnfam=\ninescaln
  \def\scaln{\fam\scalnfam\ninescaln}
  \textfont\ttfam=\ninett
  \def\tt{\fam\ttfam\ninett}%
  \textfont\bffam=\ninebf\scriptfont\bffam=\sixbf\scriptscriptfont\bffam=\fivebf
  \def\bf{\fam\bffam\ninebf}%
  \abovedisplayskip=3pt plus2pt minus 2pt
  \belowdisplayskip=3pt plus1pt minus 2pt
  \abovedisplayshortskip= 2pt plus 2pt minus 1pt
  \belowdisplayshortskip= 2pt plus 2pt minus 1pt
  \smallskipamount=2pt plus 1pt minus 2pt
  \medskipamount=3pt plus 2pt minus 2pt
  \bigskipamount=7pt plus 3pt minus 3pt
  \setbox\strutbox=\hbox{\vrule height 5pt depth 2pt width 0pt}%
  \normalbaselineskip=11pt plus.3pt minus.2pt\normalbaselines\rm}

  \def\sevenpoint{%
  \textfont0=\sevenrm\scriptfont0=\sixrm\scriptscriptfont0=\fiverm
  \def\rm{\fam0\sevenrm}%
  \textfont1=\seveni\scriptfont1=\sixi
  \scriptscriptfont1=\fivei\def\oldstyle{\fam1\seveni}%
  \textfont2=\sevensy\scriptfont2=\sixsy\scriptscriptfont2=\fivesy
  \textfont3=\sevenex\scriptfont3=\fiveex
  \textfont\itfam=\sevenit
  \def\it{\fam\itfam\sevenit}%
  \textfont\slfam=\sevensl
  \def\sl{\fam\slfam\sevensl}%
  \textfont\bbfam=\sevenbb \scriptfont\bbfam=\sixbb\scriptscriptfont\bbfam=\fivebb
  \def\bb{\fam\bbfam\sevenbb}%
  \textfont\msafam=\sevenmsa\scriptfont\msafam=\sixmsa
  \textfont\scalnfam=\sevenscaln
  \def\scaln{\fam\scalnfam\sevenscaln}
  \textfont\bffam=\sevenbf\scriptfont\bffam=\sixbf\scriptscriptfont\bffam=\fivebf
  \def\bf{\fam\bffam\sevenbf}%
  \textfont\ttfam=\seventt
  \abovedisplayskip=2pt plus2pt minus 2pt
  \belowdisplayskip=2pt plus1pt minus 2pt
  \abovedisplayshortskip= 1pt plus 2pt minus 1pt
  \belowdisplayshortskip= 1pt plus 2pt minus 1pt
  \smallskipamount=2pt plus 1pt minus 2pt
  \medskipamount=3pt plus 2pt minus 2pt
  \bigskipamount=7pt plus 3pt minus 3pt
  \setbox\strutbox=\hbox{\vrule height 5pt depth 2pt width 0pt}%
  \normalbaselineskip=9pt\normalbaselines\rm}

 \def\onzepts{%
 \textfont0=\elevenrm\scriptfont0=\ninerm
 \textfont1=\eleveni\scriptfont1=\ninei
}

\def\douzepts{%
  \textfont0=\twelverm\scriptfont0=\tenrm\def\rm{\fam0\twelverm}%
  \textfont1=\twelvei\scriptfont1=\teni%
  \textfont2=\twelvesy\scriptfont2=\tensy\scriptscriptfont2=\eightsy
  \textfont3=\twelveex
  \textfont\itfam=\twelveti
  \def\it{\fam\itfam\twelveti}%
  \textfont\bffam=\twelvebf\scriptfont\bffam=\tenbf\scriptscriptfont\bffam=\eightbf
  \def\bf{\fam\bffam\twelvebf}%
  \textfont\bbfam=\twelvebb \scriptfont\bbfam=\tenbb
  \def\bb{\fam\bbfam\twelvebb}%
  \textfont\msafam=\twelvemsa\scriptfont\msafam=\tenmsa
  \textfont\scalnfam=\twelvescaln
  \normalbaselineskip=15pt\normalbaselines\rm}

\def\quatorzepts{%
  \textfont0=\fourteenrm\scriptfont0=\twelverm\def\rm{\fam0\fourteenrm}%
  \textfont1=\fourteenib\scriptfont1=\twelveib%
  \textfont2=\fourteensy\scriptfont2=\twelvesy\scriptscriptfont2=\tensy
  \textfont3=\fourteenex
  \textfont\itfam=\fourteenti
  \def\it{\fam\itfam\fourteenti}%
  \textfont\bffam=\fourteenbf\scriptfont\bffam=\twelvebf\scriptscriptfont\bffam=\tenbf
  \def\bf{\fam\bffam\fourteenbf}%
  \textfont\bbfam=\fourteenbb \scriptfont\bbfam=\twelvebb
  \def\bb{\fam\bbfam\fourteenbb}%
  \textfont\msafam=\fourteenmsa\scriptfont\msafam=\twelvemsa
  \textfont\scalnfam=\twelvescaln
  \normalbaselineskip=18pt\normalbaselines\rm}

% Bibliographies, journaux

\def\AA{{\it Acta Arith.}}

%Insertion de figures et illustrations
\def\picture #1 by #2 (#3){\leavevmode\vbox to #2{
     \hrule width #1 height 0pt depth 0pt
      \vfill
       \special{picture #3}}}

\def\illustration #1 by #2 (#3) scaled #4{\dimen1=#2
  \divide\dimen1 by 1000
  \multiply\dimen1 by #4
  \vtop to \dimen1{\dimen1=#1
  \divide\dimen1 by 1000
  \multiply\dimen1 by #4
  \hsize=\dimen1\vss
  \noindent\special{illustration #3 scaled #4}}}

\ifx\couleurs\oui

\fi

\anglais

\optionparag=1
\def\paradouze{oui}
\vsize=255truemm
\voffset=-3truemm
\ifx\optionkeymacros\undefined\else \fi

\catcode`\Œ=\active\defŒ{{\aa}}       % option a
\catcode`\º=\active\defº{\int}        % option b (math mode) 
\catcode`\=\active\def{\c c}        % option c
\catcode`\¶=\active\def¶{\partial}    % option d (math mode)
\catcode`\Ä=\active\defÄ{\oint}       % option f (math mode) ?
\catcode`\Æ=\active\defÆ{\triangle}   % option j (math mode)
\catcode`\Â=\active\defÂ{\neg}        % option l (math mode)
\catcode`\µ=\active\defµ{\mu}         % option m (math mode)
\catcode`\¿=\active\def¿{{\o}}        % option o
\catcode`\¹=\active\def¹{\pi}         % option p (math mode w/ arg.)
\catcode`\Ï=\active\defÏ{{\oe}}       % option q 
\catcode`\§=\active\def§{{\ss}}       % option s 
\catcode`\ =\active\def {\dagger}     % option t  (math mode)
\catcode`\Ã=\active\defÃ{\sqrt}       % option v (math mode w/ arg.)
\catcode`\·=\active\def·{\Sigma}      % option w (math mode)
\catcode`\Å=\active\defÅ{\approx}     % option x (math mode)
\catcode`\½=\active\def½{\Omega}      % option z (math mode)
\catcode`\£=\active\def£{{\it\$}}     % option 3 ($ from italic font)
\catcode`\°=\active\def°{\infty}      % option 5 (math mode)
\catcode`\¤=\active\def¤{{\S}}        % option 6 
\catcode`\¦=\active\def¦{{\P}}        % option 7
\catcode`\¥=\active\def¥{\bullet}     % option 8 
\catcode`\»=\active\def»{\leavevmode\raise.585ex\hbox{\b a}}      % option 9
\catcode`\¼=\active\def¼{\leavevmode\raise.6ex\hbox{\b o}}        % option 0
\catcode`\­=\active\def­{\not=}       % option = (math mode)
\catcode`\²=\active\def²{\leq}        % option , (math mode)
\catcode`\³=\active\def³{\geq}        % option . (math mode)
\catcode`\Ö=\active\defÖ{\div}        % option / (math mode)
\catcode`\É=\active\defÉ{{\dots}}     % option ; 
\catcode`\¾=\active\def¾{{\ae}}       % option '
\catcode`\Ç=\active\defÇ{\og}         % option \ (math mode)
\catcode`\Ò=\active\defÒ{``}          % option [
\catcode`\Á=\active\defÁ{!`}          % option !
\catcode`\¢=\active\def¢{\rlap/c}     % option 4
\catcode`\Ô=\active\defÔ{`}           % option ] 
\catcode`\Õ=\active\defÕ{'}           % shift option ]

% macintosh "shift-option" generated characters

\catcode`\=\active\def{{\AA}}       % shift-option A
\catcode`\'=\active\def'{\c C}        % shift-option C
\catcode`\¯=\active\def¯{{\O}}        % shift-option O
\catcode`\¸=\active\def¸{\Pi}         % shift-option P (math mode)
\catcode`\Î=\active\defÎ{{\OE}}       % shift-option Q
\catcode`\®=\active\def®{{\AE}}       % shift-option '
\catcode`\×=\active\def×{\diamond}    % shift-option V (math mode)
\catcode`\¡=\active\def¡{\accent'27}  % shift-option 8
\catcode`\Ó=\active\defÓ{''}          % shift-option [
\catcode`\±=\active\def±{\pm}         % shift-option = (math mode)
\catcode`\È=\active\defÈ{\fg}         % shift-option \ (math mode)
\catcode`\À=\active\defÀ{?`}          % shift-option / 
\catcode`\Ð=\active\defÐ{--}          % option - (en-dash)
\catcode`\Ñ=\active\defÑ{---}         % shift-option - (em-dash)

% define the macintosh "composite" characters

\catcode`\Š=\active\defŠ{\"a}        % option u, then  a
\catcode`\'=\active\def'{\"e}        % option u, then  e
\catcode`\•=\active\def•{\"{\i}}     % option u, then  i
\catcode`\š=\active\defš{\"o}        % option u, then  o
\catcode`\Ÿ=\active\defŸ{\"u}        % option u, then  u
\catcode`\Ø=\active\defØ{\"y}        % option u, then  y
\catcode`\å=\active\defå{\^A}        %  ^, then  A
\catcode`\€=\active\def€{\"A}        % option u, then  A
\catcode`\…=\active\def…{\"O}        % option u, then  O
\catcode`\†=\active\def†{\"U}        % option u, then  U
\catcode`\‡=\active\def‡{\'a}        % option e, then  a
\catcode`\Ž=\active\defŽ{\'e}        % option e, then  e
\catcode`\'=\active\def'{\'{\i}}     % option e, then  i
\catcode`\—=\active\def—{\'o}        % option e, then  o
\catcode`\œ=\active\defœ{\'u}        % option e, then  u
\catcode`\ƒ=\active\defƒ{\'E}        % option e, then  E
\catcode`\æ=\active\defæ{\^E}        %  ^, then  E
\catcode`\é=\active\defé{\`E}        %  
\catcode`\ˆ=\active\defˆ{\`a}        % option `, then  a
\catcode`\=\active\def{\`e}        % option `, then  e
\catcode`\"=\active\def"{\`{\i}}     % option `, then  i
\catcode`\˜=\active\def˜{\`o}        % option `, then  o
\catcode`\=\active\def{\`u}        % option `, then  u
\catcode`\Ë=\active\defË{\`A}        % option `, then  A
\catcode`\‹=\active\def‹{\~a}        % option n, then  a
\catcode`\–=\active\def–{\~n}        % option n, then  n
\catcode`\›=\active\def›{\~o}        % option n, then  o
\catcode`\Ì=\active\defÌ{\~A}        % option n, then  A
\catcode`\"=\active\def"{\~N}        % option n, then  N
\catcode`\Í=\active\defÍ{\~O}        % option n, then  O
\catcode`\‰=\active\def‰{\^a}        % option i, then  a
\catcode`\=\active\def{\^e}        % option i, then  e
\catcode`\"=\active\def"{\^{\i}}     % option i, then  i
\catcode`\™=\active\def™{\^o}        % option i, then  o
\catcode`\ž=\active\defž{\^u}        % option i, then  u

\let\optionkeymacros\null

\def\gS{{\got S}}
\font\tengp=cmmib10

\font\ninegp=cmmib10 at 9pt
               \font\eightgp=cmmib8
              \font\sevengp=cmmib7
%  \font\sixgp=cmmib6
               
               \newfam\gpfam
               \textfont\gpfam=\tengp\scriptfont\gpfam=\sevengp
               \def\gp{\fam\gpfam\textfont1=\tengp\scriptfont1=\sevengp\tengp}%

\font\tenib=cmmib10
         \font\nineib=cmmib10 scaled 900
\font\sevenib=cmmib10 scaled 700
         \font\fiveib=cmmib10 scaled 500
\def\itg{\ib}
\def\0{{\bf 0}}
\def\a{{\itg a}}
\def\b{{\itg b}}
\def\cc{{\itg c}}
\def\m{{\itg m}}
\def\n{{\itg n}}
\def\s{{\itg s}}
\def\x{{\itg x}}
\def\y{{\itg y}}
\def\R{{\itg R}}
\def\el{\hbox{$\gp \ell$}}

\def\Ns{{\N^*}}

\def\dens{\mathop{\rm dens}\nolimits}

\def\lcm{\mathop{\rm lcm}\nolimits}
\def\gcd{\mathop{\rm gcd}\nolimits}
\def\quadr#1#2#3#4{{{\sct#1\atop\sct#2}\atop\sct#3}\atop\sct#4}

\def\gs{{\got s}}

\def\sgxi{\hbox{$\sct\gp\xi$}}

\font\tenib=cmmib10
             \font\nineib=cmmib10 scaled 900
\font\sevenib=cmmib10 scaled 700
             \font\fiveib=cmmib10 scaled 500
\def\itg{\ib}

\def\auteur{R. de la Bretche \& G. Tenenbaum}
 
\def\titrart{Mean values of arithmetic functions and application to sums of powers} 

\hautspages{\auteur}{\titrart}
%\dateheure
\titrecentre{\titrart}
\bigskip\medskip
\centerline {\auteur} 
\bigskip\bigskip
{\eightpoint\leftskip1cm\rightskip1cm
\noi{\bf Abstract.} We provide new upper bounds for sums of certain arithmetic functions in many variables at polynomial arguments and, exploiting recent progress on the mean-value of the Erd\H os-Hooley $\Delta$-function, we derive lower bounds for the cardinality of those integers not exceeding a given limit that are expressible as  certain  sums of powers.
 \PAR
\medskip
\noi
{\bf Keywords:} Erd\H os-Hooley $\Delta$-function, sieve, sums of powers, arithmetic functions of many variables, mean-values of arithmetic functions.\par
\smallskip 
\noi \bf 2020 Mathematics Subject Classification: \rm primary   11N25; secondary 11N37, 11N64, 11P05.\par }
\bigskip
\medskip
\paraun{Introduction  } 
The  Erd\H os-Hooley $\Delta$-function is defined by the formulae
 $$\Delta(n,u):=|\{d:d\mid n,\,\e^u<d\leqslant \e^{u+1}\}|,\quad\Delta(n):=\max_{u\in\r}\Delta(n,u)\qquad (n\in \N^*).$$
 Recent, spectacular progress on the mean-value of the $\Delta$-function has been obtained in \citer{KT23}, \citer{FKT23}, and refined in  \citer{dlBT24}. This  invites one to revisit some of the applications linked to this quantity. Put
 $$S(x):=\sum_{n\leqslant x}\Delta(n),\quad \gS(x):={1\over \log x}\sum_{n\leqslant x}{\Delta(n)\over n},$$
 so that (see   \citeplus{HT88}{th. 61})
 $$S(x)\ll x\gS(x)\qquad (x\geqslant 2).\note{ It is actually easy to show that $S(x)\asymp x\gS(x)\ (x\geqslant 2)$, but we shall not need the lower bound here.}$$
 After the works of Erd\H os \citer{Er73}, Hall-Tenenbaum \citer{HT82}, \citer{HT88}, Hooley \citer{Ho79}, Tenenbaum \citer{Te85}, Koukou\-lopoulos-Tao \citer{KT23}, Ford-Koukoulopoulos-Tao \citer{FKT23}, and La Bretche-Tenenbaum \citer{dlBT24}, we are now equipped with the bounds
  $$(\log_2x)^{1+\eta+o(1)}\ll \gS(x)\ll (\log_2x)^{5/2}\qquad (x\geqslant 3),\note{Here and throughout, we let $\log_k$ denote the $k$th iterated of the logarithm.}\eqdef{estBT}$$ 
  where $\eta\approx0.353327$ is the exponent appearing in the new lower bound for the normal order of $\Delta(n)$ by Ford, Green and Koukoulopoulos \citer{FGK23}. 
 A further description of successive advances is available in the article \citer{dlBT24}.
\par \medskip
The fact that the asymptotic order of magnitude of $\gS(x)$ is a power of $\log_2x$ makes it pertinent to refine the estimates relating to the various applications, as the number of representations of an integer as a sum of powers, the Diophantine approximation of an irrational number by a rational with square denominator, or Waring's problem. \par 
Two other applications are also worth mentioning. First note that \citeplus{NT98}{cor. 8} readily yields that, for any given $\varepsilon>0$, and uniformly for $x\geqslant 3$, $1\leqslant |a|\leqslant x$, $x^\varepsilon<y\leqslant x$, we have
$$\sum_{x<p\leqslant x+y}\Delta(|p+a|)\ll{|a|y\gS(x)\over  \varphi(|a|) \log x}\ll{|a|y(\log_2x)^{5/2}\over \varphi(|a|)\log x},$$
where the letter $p$ denotes a prime number.
Second, observe  that the upper bound of \eqref{estBT} furnishes, for any fixed $\varepsilon>0$ and uniformly for $ 1\leqslant a \leqslant q\leqslant x^{2-\varepsilon}$, $(a,q)=1$, $x\geqslant 3$, the estimate
$$\sum_{\di{n,m\leqslant x}{mn\equiv a \mod q}} 1\ll {x^2\over q}(\log_2x)^{5/2}.\eqdef{inegDeshou}$$
Indeed, if $N\leqslant x^{1-\varepsilon/2}$, a trivial bound, using $\tau(n)=n^{o(1)}$ suffices, and if $x^{1-\varepsilon/2}<N\leqslant x$, we have
$$\sum_{\tri{ m\leqslant x}{N<n\leqslant 2N}{mn\equiv a \mod q}} 1\leqslant \sum_{0\leqslant k\leqslant 2Nx/q} \Delta(kq+a)
\ll {Nx\over q} \gS(x) ,$$
where the latter bound results from a special case of \citeplus{NT98}{th. 1}. This implies \eqref{inegDeshou} by summing over $N$. 
\goodbreak
We intend here to return to the problem of  sums of powers, also addressed in \citer{KT23}.
\par 
Given $\cc:=\{c_j \}_{0\leqslant j\leqslant t},\,\el:=\{\ell_j \}_{0\leqslant j\leqslant t}\in(\N^*)^{t+1}$,  consider the number $$r(n) = r(n;\cc,\el) $$
of representations of a natural integer $n$ in the form 
$$n=
\sum_{0\leqslant j\leqslant t} c_jm_j^{\ell_j} \qquad (\m \in \N^{t+1}).$$
 In all the sequel, we shall assume without loss of generality that $\{\ell_j \}_{0\leqslant j\leqslant t}$ is non-decreasing.
  \par \goodbreak
Under the condition $\sum_{0\leqslant j\leqslant t}1/\ell_j=1$, we trivially have  
 $$\qquad V_0(x;\cc,\el):=\{ n\leqslant x : r(n)\geqslant 1\}\leqslant x\qquad\qquad(x\geqslant 1).$$
The exact order of magnitude of $V_0(x)=V_0(x;\cc,\el)$ is only known in the case $\el=(2,2)$.  A famous conjecture (see, for instance, \citer{BB95}, or \citer{W18}) is that $V_0(x)\asymp x$ for $\el=(3,3,3)$, $\cc=(1,1,1)$. 
 The case $t=2$, $\cc=(1,1,1)$, $\el=(2,4,4) $ was studied by Hooley \citer{Ho79}, Tenenbaum  \citer{Te86}, Robert \citer{R11}, and Koukoulopoulos--Tao~\citer{KT23}.  Robert \citer{R11} also considered the case  $\el=(2,3,6)$. 
  \par  
We obtain the following result.
\Propt{thRobert}{Let $t  \geqslant 2
$, $\cc\in(\N^*)^{t+1}$,  $\el=\{\ell_j \}_{0\leqslant j\leqslant t}$ with $\ell_0=2$, $\ell_1\in \{ 3,4\}$ and $\sum_{j=1}^t 1/\ell_j=\dm$. Let $s:=\max\{j\in[1,t-1]:\ell_{t-j+1}=\ldots=\ell_t\}$. Assume:\par 
{\rm(i)} $1\leqslant s\leqslant 3$, 
\par 
{\rm(ii)} $\ell_t\geqslant 26$ if $s=3$, \par 
{\rm(iii)} $\sum_{j\geqslant r} 1/\ell_j\leqslant 1/\ell_{r-1}\ (1\leqslant r\leqslant t-s+1).$ 
\par 
Then, we have 
$$V_0(x;\cc,\el)\gg x/\gS(x)\gg x/(\log_2x)^{5/2}.\eqdef{estV0}$$
 }
 We have not sought to specify uniformity in respect to the coefficients $c_j$ in this statement.
 \par 
 Observe that the following sets fulfill the assumptions:\par \smallskip
{\hglue5mm $\eqalign{&t=2,\,\el\in \{  (2,3,6) ,\, (2,4,4) \}, \cr
 &t=3,\,\el\in \{  (2,3,7,42),\, (2,3,8,24) , \,  (2,3,9,18), \,  (2,3,10,15),\,  (2,3,12,12),\, (2,4,8,8)   \},\cr
  &t=4,\,\el\in\{ (2,3,7,84,84), \, (2,3,8,48, 48),% \,  (2,3,18,18,18), 
   \, (2,4,8,16,16),\,  (2,3,12,24,24)\},  \cr
     &t=5,\,\el\in\{(2,3,12,36,36,36),\,  (2,3,10,45,45,45),\, % (2,4,8,24,24,24),\,  
     (2,4,8,16,32,32)\},}$}
\par \smallskip
\noi
this list being non exhaustive. New sets can indeed be constructed by increasing the cardinality~$t$. For instance, when  $\el$ belongs to the above list, we can define $\el^+\in(\N^*)^{t+2}$ by setting $\ell^+_j:=\ell_j$ for $0\leqslant j\leqslant t-1$, $\ell_t^+:=2\ell_t$, $\ell_{t+1}^+=2\ell_t$. Yet another possibility is to define $\el^*\in(\N^*)^{t+3}$ with $\ell^*_t=\ell_{t+1}^*=\ell_{t+2}^*:=3\ell_t$ provided $\ell_t\geqslant  9$. Thus, $\ell=(2,3,6)$ induces $ (2,3,12,12) $
and  $ (2,3,18,18, 18). $
\par 
 The article \citer{KT23} comes back on the case $\el=(2,4,4)$---see, however, the remark following the statement of \ref{majV2} {\it infra}. \par 
 The cases with $t\geqslant 3$ in \ref{thRobert} are new. The corresponding results rely on the bound
  $$\Vbs{\bigg\{(m_{t-2}, m_{t-1}, m_{t},n_{t-2}, n_{t-1}, n_{t}):\sum_{t-2\leqslant j\leqslant t} c_jm_j^{\ell_t} =\sum_{t-2\leqslant j\leqslant t} c_jn_j^{\ell_t}\leq x\bigg\}}\ll x^{3/\ell_t},\eqdef{Salb}$$
  established by 
 Salberger provided $\ell_t\geq 26$---see the proof of \ref{majV2=} {\it infra} for precise references. In   private communication \citer{Sa24}, Salberger informed us that, as a direct consequence of his works \citer{Sa23} and \citer{Sa05}, the bound \eqref{Salb} persists under the weaker condition $\ell_t\geq 16$. Assuming this, we get that \eqref{estV0} also holds for $\el=(2,3,18,18,18)$ and $\el=(2,4,8,24,24,24)$.  
\medskip
 The starting point of the proof of \ref{thRobert} Êis classically the  Cauchy-Schwarz inequality
 $$V_1(x)^2\leqslant V_0(x)V_2(x)\eqdef{inegCSV}$$
where we have set
 $$V_j(x)=V_j(x;\cc,\el):=\sum_{n\leqslant x} r(n)^j\qquad (j=1,2).\eqdef{defVj}$$
We shall deduce \ref{thRobert} from a general estimate regarding $j=2$.
\par 
\paraun{Notation}
 Some notation must be introduced in order to state our results.  \par
Given $k\in\N^*$, $\a=(a_j)_{1\leqslant j\leqslant k},\b=(b_j)_{1\leqslant j\leqslant k}\in\N^{*k}$, we put $$\gs\a:=\sum_{1\leqslant j\leqslant k}a_j,\qquad \wp\a:=a_1\cdots a_k,\quad\a\b=(a_1b_1,\ldots,a_kb_k).$$ \par 
For each value of the parameters $A\geqslant 1$, $B\geqslant 1$,
$\varepsilon>0$, we designate by
$\M_k(A,B,\varepsilon)$ the class of those non-negative arithmetic functions
 $F:\N^{*k}:\to\r^+$,  satisfying the condition
$$F(\a\b)\leqslant \min
\big\{A^{\Omega(\wp\a)}, B(\wp\a)^\varepsilon \big\}F(\b)\eqdef{defMk}$$ 
for all $\a,\b\in\N^{*k}$ such that
$(\wp\a, \wp\b)=1$. By convention, the functions $F\in\M_k(A,B,\varepsilon)$ are extended to $\N^k$ by setting $F(\a)=0$ if $\min_ja_j=0$.
When $F\neq0$, we define
$$G(\a)=G_F(\a):=\max_{\tri{\b\in\N^{*k}}{(\wp\a,\wp\b)=1}{F(\b)\neq
0}}{F(\a\b)\over
F(\b)}\qquad
\big(\a\in\N^{*k}\big).\eqdef{defG}$$
\par
Given a family $\{Q_j\}_{j=1}^{k}\in
\z[X_1, \ldots, X_t]^k$  of polynomials in $t$ variables, we put
$$Q=\prod_{1\leqslant j\leqslant
k}Q_j=\prod_{1\leqslant h\leqslant
r}R_h^{\gamma_h}\in\z[X_1, \ldots, X_t], \quad g:=\deg
Q,\eqdef{decompo}$$
where the $R_h$ are irreducible in $\z[X_1, \ldots, X_t]$.
We may then write, canonically,
$$Q_j=\prod_{1\leqslant h\leqslant
r}R_h^{\gamma_{jh}}\qquad (1\leqslant j\leqslant
k)$$
where $\gamma_{jh}\geqslant 0$ for all $j,h$, so that
$$\gamma_h=\sum_{1\leqslant j\leqslant
k}\gamma_{jh}\qquad (1\leqslant h\leqslant r), \qquad \sum_{1\leqslant h\leqslant r} \gamma_h \deg R_h =g.$$
\par
Recall that a polynomial is said to be primitive if the greatest common divisor of its coefficients is 1. In all the sequel we assume the form $Q$ to be primitive, which implies that the same holds for all the $Q_j$.
\par \goodbreak
For $T\in\z[X_1, \ldots, X_t]$, we write
$$
\varrho_T^+(s):=\sum_{\di{\sgxi\in[1,s]^t}{T(\sgxi)\md0 s}}1\qquad
(s\geqslant 1).\eqdef{defrho+}$$
Letting $\kappa(s):=\prod_{p|s}p$ denote the squarefree kernel
of a natural integer $s$, we further put
$$\K(\s):=\lcm(s_j\kappa(s_j))_{1\leqslant j\leqslant r}=[s_1\kappa(s_1),\ldots,
s_r\kappa(s_r)]\qquad
\big(\s=(s_1,\ldots,s_r)\in\N^{*r}\big).\eqdef{defKn}$$
 Given $\R\in\z[X_1,\ldots,X_t]^r$,  we may then define two arithmetic functions of $\s=(s_1,\ldots,s_r)\in\N^{*r}$~by $$\leqalignno{\rho^+_\R(\s)&:=\sum_{\di{\sgxi\in [1,\wp\s]^t}{R_h(\sgxi)\md0{s_h}\ (1\leqslant
h\leqslant r)}}1,&\eqdef{defrho++}\cr
\varrho^\#_\R(\s)&:=\sum_{\tri{\sgxi\in[1,\K(\s)]^t }{s_h\|R_h(\sgxi)\
(1\leqslant h\leqslant r)}{(R_h(\sgxi)/s_h,\wp \s)=1   } }1,&\eqdef{defrhorho}\cr}
$$  
where the symbol $a\| b$ means that the conditions
 $a|b$ and $(a,b/a)=1$ hold simultaneously. It is worthwhile to note that the above definition of $\varrho^\#_\R$ is more restrictive than that of \citer{dlBT12}. The grounds for this alteration are described in \citer{H14}.\par
For any
$\s=(s_1,\ldots,s_r)\in\N^{*r}$, let us put
$$s'_j:=\prod_{1\leqslant h\leqslant
r}s_h^{\gamma_{jh}}\quad (1\le j\le k),\quad\s':=(s'_1,\ldots,s_k'),\qquad
s'':=\wp\s'=\prod_{1\leqslant j\leqslant
k}s'_j=\prod_{1\leqslant h\leqslant %R n->s
r}s_h^{\gamma_h}.
$$
Thus, given $\m=(m_1,\ldots, m_t)\in\N^{*t}$, and writing
  $s_h=R_h(\m)\ (1\le h\le r)$, we have
$$s'_j=Q_j(\m)\quad(1\leqslant j\leqslant k),\quad s''=Q(\m).$$
With these pieces of notation, to any function
 $F$ in $\M_k(A,B,\varepsilon)$ we may  associate a further
function $\wh F$ in $\M_r(A^g, B,g\varepsilon)$
defined by
$$\wh F(\s)=F(\s').$$
We then designate by  $\wh G:=G_{\wh F}$ the function associated to $\wh F$ by \eqref{defG} with $k=r$ and, for $1\leqslant h\leqslant r$,
we let  $\wh G_h:\N^*\to\r^+$ denote the composition of $\wh
G$ with the $h$th  coordinate.
\medskip
 \paraun{Statement of results}
The key tool of this work consists in an estimate of the quantityÊ
$$\eqalign{S&:=\sum_{x_j<n_j\leqslant x_j+y_j \, (1\leqslant j\leqslant t)}
F\big(|Q_1(\n)|,\ldots, |Q_k(\n)|\big)
=\sum_{x_j<n_j\leqslant x_j+y_j \, (1\leqslant j\leqslant t)}
\wh F\big(R_1(\n),\ldots, R_r(\n)\big).\cr}
$$
The following result is obtained. It is the $t$-dimensional counterpart of that of \citer{H12}. We let $\|Q\|$ denote the maximum of the absolute values of the coefficients of a real form $Q$.
\par 
\Propt{th}{Let $k,t\in\Ns$, and let  $\{Q_j
\}_{j=1}^{k}\in \z[X_1, \ldots, X_t]^k$ be a family of primitive polynomials. Define
$Q\in\z[X_1, \ldots, X_t]$, $\{R_h\}_{1\leqslant h\leqslant r}\in\z[X_1, \ldots, X_t]^r$,
$g\in\N$ by \eqref{decompo}, $\varrho_Q^+$ by
\eqref{defrho+}, and $\varrho_\R^\#$ by~\eqref{defrhorho}. For all
        $$\alpha\in]0,1],\quad \beta\in ]0,1[,\quad
A\geqslant 1,\quad B\geqslant 1, \quad 0<\varepsilon\leq
\alpha \beta/ \{50g^2( \beta g+1)\},\eqdef{pourtous}$$ and uniformly under the
 conditions
 $$\eqalign{&F\in  \M_k(A,B,\varepsilon),\quad\x=(x_1,\ldots,x_t)\in\N^{*t},\quad\y=(y_1,\ldots,y_t)\in\N^{*t},\cr
&x:=\min_{1\leqslant j\leqslant t} x_j\geqslant c\big\{\max_{1\leqslant j\leqslant t}x_j+\no{Q}\big\}^{ \beta},\quad
 x_j^\alpha\leqslant y_j\leqslant
x_j \quad(1\leqslant j\leqslant t),\cr}\eqdef{condth}$$we have
$$ \eqalign{\sum_{\di{\n\in\N^{*t}}{x_j-y_j<n_j\leqslant x_j} (1\leqslant j\leqslant t)}&
F\big(|Q_1(\n)|,\ldots, |Q_k(\n)|\big) \ll
 \wp\y  E_\R(\gs \x)\prod_{g< p\leqslant x}\Big(1-{\varrho_Q^+(p)\over
p^t}\Big)\cr}\eqdef{estfond}$$
where  $c$ and the implicit constant depend at most upon $g$, $\alpha$, $ \beta$, $A$, $B$ and where
we have set 
$$ E_\R(v):=\sum_{\di{\s\in \Ns^r}{\wp\s\leqslant v}}\wh
F(\s){\varrho_\R^\#(\s)\over
\K(\s)^t}\qquad (v\geqslant 1).\eqdef{defER}$$
}

 \rems (i) We have
 $${ \rho^\#_\R(\s)\over \K(\s)^t}\leqslant {{\varrho^+_\R}(\s)\over  (\wp\s)^t}\quad\qquad
\big(\s=(s_1,\ldots,s_r)\in{\N^{*r}}\big),\eqdef{majdiesepar+}$$
however the sum $ E_\R(x)$  may turn out to be smaller than that which follows from this upper bound. For instance, if $k=1$, $t=1$, $Q(X)=Q_1(X)=X(X+\ell)$, we have $\rho^\#_\R(p)=0$ as soon as  $p\mid \ell$, while $\rho^+_\R(p)=1$. As noted in \citer{H12} and \citer{H14}, this can be significant when it comes to obtaining  results uniform with respect to $\ell$.
 
 (ii) A weaker version of \eqref{estfond}  for $k=1$  has recently been established \citeplus{CKPS24}{th. 1.15}, where,  in definition \eqref{defER}, the left-hand side of \eqref{majdiesepar+} is replaced by its majorant. Note that the dependence of our bounds with respect to the coefficients of $Q$ is  effective.
  
 (iii) In the case $t=1$,  \ref{th} is due to Henriot \citer{H12}, \citer{H14}, extending \citer{Sh80},  then  \citer{Na92} and \citer{NT98}. When the $Q_j$ are binary forms in $\z[X_1,X_2]$,  \ref{th} has been proved in \citer{dlBT12} with a slightly different definition of $\rho^\#$, refining and extending the estimates of  \citer{dlBB06}.
 
 (iv) When $F$ is multiplicative and satisfies an assumption such that $F(\n)\geqslant \eta^{\Omega(\wp\n)}$ for some $\eta>0$, it is possible to establish lower bounds for the the left-hand side of   \eqref{estfond} 
 of  the same order of magnitude as the right-hand side. We refer to \citer{H12} for this aspect of the question. 
 
 \paraun{Proof of \ref{th}}

Recall a consequence of the Schwartz-Zippel lemma --- see \citeplus{S80}{lemma 1}, or \citeplus{Z79}{th. 1}. 
\Propl{ZS}{Let $t\geqslant 1$ and let $Q\in \z[X_1, \ldots, X_t] $ be a primitive polynomial of degree $g$.  For all prime numbers $p$, we have
$$\rho_Q^+(p )\leqslant g p^{ t-1 }.\eqdef{majZS}$$In particular, for all integers $\nu\geqslant 1$,  we have $\rho_Q^+(p^\nu)\leqslant g p^{\nu t-1}<p^{\nu t}$ provided $p>g$.}
\medskip \goodbreak  
  Let $c(Q)$ denote the content of the polynomial $Q$ of $t$ variables and let $v_p(n)$ denote the $p$-adic valuation of an integer $n$.  %  \Propl{lemme1}{Let $Q\in \z[X] $ a polynomial of degree $g$.  For all prime numbers $p$ and all integers $\nu\geqslant 1$, we have $$\rho_Q^+(p^\nu)\leqslant g p^{\nu (1-1/g)}p^{v_p(c(Q))/g}.\eqdef{majlemme1}$$}
   The following lemma is an effective version of \citeplus{CKPS24}{lemma 2.8}---see also \citeplus{PSW16}{lemma 4.10}.
   
  \Propl{lemme1}{Let $t\geqslant 1$ and let $Q\in \z[X_1, \ldots, X_t] $ be a   polynomial of degree $g$.  For all prime numbers $p$ and all integers $\nu\geqslant 1$, we have
$$\rho_Q^+(p^\nu)\leqslant g^t(\nu+1)^{t-1} p^{\nu (t-1/g)+ v_p(c(Q))/g}.\eqdef{majlemme2}$$ }

   \nid Consider first the case $t=1$. When $(c(Q),p)=1$,  inequality \eqref{majlemme2} has been proved by Stewart \citer{S91}. If $\gamma=v_p(c(Q)) \geqslant 1 $, we have
$$\rho_Q^+(p^\nu)=p^\gamma \rho_{Q^*}^+(p^{\nu-\gamma})  \leqslant g p^{\nu (1-1/g)}p^{v_p(c(Q))/g}.$$
with $Q^*:=Q/p^\gamma$.  Thus \eqref{majlemme2} holds for $t=1$. 
 
   We now proceed by induction on $t$. Let $x_1$ be such that $Q$ depends on $x_1$ and write 
   $$Q(x_1,x_2,\ldots, x_t)=\sum_{0\leqslant j\leqslant g} Q_j( x_2,\ldots, x_t)x_1^j $$
   with $Q_j \in \z[X_2, \ldots, X_t] $  and $g_j:=\deg Q_j\leqslant g-j$.
   By  the case $t=1$, we have 
       $$\rho_Q^+(p^\nu)\leqslant g p^{\nu (1-1/g)}\sum_{0\leqslant \gamma\leqslant \nu}
       p^{\gamma/g+(t-1)(\nu-\gamma)} N(\gamma,p)$$
       with $N(\gamma,p):=|\{ \x\in ( \z/p^\gamma\z)^{t-1}: Q_j(\x)\md0{p^\gamma}\ (0\leqslant j\leqslant g)\}|$.
  \par  Now the induction hypothesis yields
 $$\eqalign{ N(\gamma,p)
& \leqslant \min_{0\leqslant j\leqslant g}
|\{ \x\in ( \z/p^\gamma\z)^{t-1} :  Q_j(\x)\md0{p^\gamma}  \}|
 \cr& \leqslant \min_{0\leqslant j\leqslant g} g_j^{t-1}(\nu+1)^{t-2} p^{\gamma (t-1)+\min\{ v_p(c(Q_j))-\gamma, 0\}/g_j}
  \cr& \leqslant   g^{t-1}(\nu+1)^{t-2} p^{\gamma (t-1-1/g )+v_p(c(Q))/g },}$$ since 
  $ \min_{0\leqslant j\leqslant g} v_p(c(Q_j))=v_p(c(Q))$.   This completes the proof of  \eqref{majlemme2}.   \qed
 \bigskip
 \goodbreak
%  The conditions appearing \eqref{densrhod} may also be rewriten as 
% $$a_h=a_h'a_h'',\quad a_h'\mid R_h(n),\quad (a_h,
%R_h(n)/a_h')=1,\quad  (p\mid \wp\a, \, p\nmid a_h)\Rightarrow p\nmid R_h(n)  \quad (1\leqslant h\leqslant
%r).$$
\medskip
Let us now state the analogue in our context of \citeplus{H12}{lemma 6}. This is a straightforward consequence of the fundamental lemma of the combinatorial sieve. It coincides with \citeplus{dlBT12}{lemma 3.3} when $t=2$  and with \citeplus{H12}{lemma 6} when $t=1$.   As in  \citeplus{H12}{lemma 6}, we let $\alpha$ be an arbitrary constant in $]0,1]$ and define
$$ {{\alpha_1:=\ft3{25} \alpha,\quad% \varepsilon_2=\ft13\varepsilon_1,\quad
 \alpha_2:=\varepsilon_1/6g}}.\eqdef{defeps}$$

\Propl{BT6}{Let $t\geqslant 1$, $0<\alpha<1$, and let $Q, R_h\in\z[X_1, \ldots, X_t]$ be polynomials satisfying the hypotheses of \ref{th}. Under the conditions 
 $$\a\in\Ns^r,\ \x,\,\y\in\N^{*t},\ w:=\min_jx_j,\  \wp\a\leq
w^{\alpha_1},\  z\leqslant w^{\alpha_2},\ x_j^{\alpha}\leq
y_j\leqslant x_j,$$ we have
$$\sum_{\quadr{x_j<n_j\leqslant x_j+y_j \, (1\leqslant j\leqslant t)}{a_h\| R_h(\n)\ (1\leqslant
h\leqslant r)}{(p\mid Q(\n), \, g< p\leqslant z
)\Rightarrow p\mid \wp\a}
{(R_h(\n)/a_h,\wp\a)=1\,(1\leqslant h\leqslant r)}}1\ll
\wp\y{\varrho_\R^\#(\a)\over \K(\a)^t}\prod_{\di{g<
p\leqslant z}{p\,\nmid\, \wp\a}}\Big(
1-{\varrho_Q^+(p)\over p^r}\Big).\eqdef{majBT6}
$$}
Note that, since the product in \eqref{estfond} only covers primes $p>g$, it  avoids any fixed prime divisor of the~$R_h$. It is also worthwhile to keep in mind that
$${\rho_\R^\#(\a)\over
\K(\a)^t}=\dens \Big\{ \n\in \z^t: a_h\parallel R_h(\n)\  (1\leqslant h\leqslant r), \, (p\mid \wp\a, \, p\nmid a_h)\Rightarrow (p\nmid R_h(\n))   \Big\} .\eqdef{densrhod}$$

 \par\goodbreak 
 \nid In the above statement, the parameters $\alpha_1$ and $\alpha_2$  defined in \eqref{defeps} coincide respectively with $\varepsilon_1$ and $\varepsilon_3$ from \citeplus{H12}{lemma~6}.
 \ref{BT6} follows  by  applying the combinatorial sieve to the sequence 
 $${\scal A} :=\left\{ Q(n_1,\ldots, n_t)\,: \eqalign{&x_j<n_j\leqslant x_j+y_j \ (1\leqslant j\leqslant t),\  a_h\| R_h(\n)\ (1\leqslant
h\leqslant r),\cr&
(p\mid \wp\a/ a_h, \, p\,\nmid\,a_h)\Rightarrow (p\,\nmid\, R_h(n))}
\right\}.$$
Thus, we can select $A_1=g+1$ in the hypotheses of the sieve---see \citeplus{HR74}{lemma~2.2},  condition~$(\Omega_1)$---
since the upper bound \eqref{majZS} implies that 
$${\rho_Q^+(p)\over p^t}\leqslant {g \over p}\leqslant 1-{1\over g+1}\quad\big(p\geqslant g+1,\,p\nmid \wp\a\big).$$ 
Further details are omitted.
\qed
\medskip
 
 Substituting \ref{BT6} to \citeplus{H12}{lemma 6} in the proof of \citeplus{H12}{th.\thinspace5} readily furnishes  estimate~\eqref{estfond}.
Note that \ref{lemme1} is crucial in order to get effective bounds in terms of $\no Q$.

\medskip 
 \paraun{Proof of \ref{thRobert}} 

\paradeuxb{Statements}
It follows from \eqref{defVj} that
$$V_2(x;\cc,\,\el)=\card\Big\{ (\m,\n)\in (\N^*)^{t+1 }\times (\N^*)^{t+1 }: \sum_{0\leqslant j\leqslant t}c_jm_j^{\ell_j}= \sum_{0\leqslant j\leqslant t}c_j n_j^{\ell_j}\leqslant x\Big\},
$$
and so  $V_2(x)=V_2(x;\cc,\,\el)=V_2^{\neq}(x)+V_2^{=}(x)$, where the subsums  $V_2^{\neq}(x)$ and $V_2^{=}(x)$ respectively correspond to the extra conditions $m_0\neq n_0$ and $m_0=n_0$. We obtain the following estimates.
\par 
\Propp{majV2}{Let $t\geqslant 1$, $\delta>0$, $\cc,\,\el\in(\N^*)^{t+1}$, satisfy $\ell_0\geq2$ and  $\delta=\sum_{1\leqslant j\leqslant t} 1/\ell_j $.
Then 
$$ V^{\neq}_2(x)\ll x^{2\delta}\gS(x)\ll x^{2\delta}(\log_2x)^{5/2} \qquad (x\geqslant 3).\eqdef{estV2neq}$$
 }
The second bound follows from \eqref{estBT}.
The exponent $5/2$ hence measures, in the current state of knowledge, the deviation of this general estimate from optimality.
 In \citeplus{KT23}{\S\thinspace8}, Koukoulopoulos and Tao indicate the exponent $=2^{4L}+15/4$ for $L=\max\ell_j$, $\ell_0=2$, $\delta=\dm$. The condition $m_0\neq n_0$ has been omitted in \citer{KT23}. Since it necessary for a general result,  the argument given in \citer{KT23} only yields  bound \citeplus{KT23}{(1.5)} in the historical case $\el=(2,4,4)$ considered by Hooley. The exponent stated in~\citer{KT23} is then $65539+3/4.$
\par 
\smallskip 
 Our approach to \ref{thRobert} is based on an inductive argument described in the following statement. 
\Propp{majV2=}{Let $t\geqslant 2$, $\cc {=\{c_j\}_{j=0}^{t}\in(\N^*)^{t+1},\el=\{\ell_j\}_{j=0}^{t}\in(\N^*)^{t+1}}$ with $\ell_0\geqslant 2$.  
 Put $\check{\cc}:=(c_1,\ldots, c_t)$, $\check{\el}=(\ell_1,\ldots, \ell_t)$. 
We have
$$V^{=}_2(x;\cc,\,\el)\leqslant \Big({x\over c_0}\Big)^{1/\ell_0} V_2(x;\check{\cc},\check{\el}) \qquad (x\geqslant 3).\eqdef{estV2=}$$
Moreover, if $t\in\{1,2\}$ and $\ell_{t-1}=\ell_t\geqslant 3$,  or if $ t=3$ and $\ell_{t-2}=\ell_{t-1}=\ell_t\geqslant 26$, then  
$$V_2(x;\check{\cc},\check{\el})\ll x^{t/\ell_t }.\eqdef{estV2check}$$}

\nid   The  bound \eqref{estV2=} is clear, and so is \eqref{estV2check} for $t=1$. 
 The remaining  bounds \eqref{estV2check} are due to Salberger: for $t=2$, see \citeplus{Sa23}{cor.\thinspace0.7}; for $t=3$, $\ell_t\geqslant 33$, see \citer{BHB04}; for $t=3$,  $\ell_t\geqslant 26 $, see  \citeplus{Sa05}{cor. 4.6}.
 %;  for $t=3$, $\ell_t\geqslant 16$, see \citer{Sa24}, explaining that  the required estimate readily follows from \citer{Sa23} and \citer{Sa05}.   
  \qed
\medskip
 As mentioned in the introduction, we have been informed by Salberger that he intends to write up the proof that \eqref{estV2=} holds for $t=3$, $\ell_t\geqslant 16$, and hopefully $\ell_t\geqslant 12$. 
\medskip
 \paradeuxb{A lemma from the literature} 

Define $P^+(n)$ as the largest prime factor of an integer $n$, with the convention that $P^+(1):=1$. The following statement is established in  \citeplus{H12}{lemmas 2 \& 5}.
\Propl{lemma2}{Let $r, g\geqslant 1$,  $A\geqslant 1$, $B\geqslant 1$,
$\varepsilon>0$  and let $\{Ê\sigma_j\}_{j=1}^r$, $\{Ê\vartheta_j\}_{j=1}^r$ be two families of multiplicative arithmetical functions.  Assume that, uniformly for prime powers $p^\nu$, we have $\sigma_j(p^\nu)\ll 1$,   $\vartheta_j(p^\nu )=1+O(1/p)$. Let  $\wh F\in\M_r(A^g, B,g\varepsilon)$, and define
$$H(\s) :=\wh F(\s) \prod_{1\leqslant j\leqslant r} \sigma_j(s_j){\varrho_\R^\#(\s)\over
\K(\s)^t}\qquad (\s\in \N^{*r}).$$
Then, for $x\geqslant 3$, we have
$$\leqalignno{ \sum_{\di{ \s\in  (\N^{*r}}{ \wp\s\leq x}} H(\s)\prod_{1\leqslant j\leqslant r} \vartheta_j(s_j)&\ll  % 
\sum_{\di{ \s\in  (\N^{*r}}{ \wp\s\leq x}} 
H(\s) ,&\eqdef{maj1sumH}
\cr
\sum_{\di{\s\in (\N^{*r}}{P^+(\wp\s)\leq x}} H(\s) &\ll  \sum_{\di{\s\in(\N^{*r}}{\wp\s\leq x}} H(\s) 
.&\eqdef{maj2sumH}}$$
}\medskip
\paradeuxb{Proof of \ref{majV2}.}

Consider $$Q(X_1,\ldots, X_t, Y_1,\ldots, Y_t):=\sum_{1\leqslant j\leqslant t} c_jX_j^{\ell_j}-\sum_{1\leqslant j\leqslant t} c_jY_j^{\ell_j} .\eqdef{nouvQ}$$
 The following statement gathers the necessary estimates regarding  function $\varrho_Q$ defined in \eqref{defrho+}.
\par \Propl{lemmeRobert}{Let $t\geqslant 1$, $\cc,\,\el\in\N^{* t}$ be such that $$
L:=\max_{1\leqslant j\leqslant t}\ell_j,\quad\gcd(c_1,\ldots,c_t)=1.$$ 
For all prime numbers $p$ and all integers $\nu\geqslant 1$ and $Q$ defined by \eqref{nouvQ}, we have
 $$\rho_Q(p )=p^{2t-1}+O(p^t), \qquad \rho_Q(p^\nu)\leqslant 
L p^{ 2t\nu-1 } , \qquad \rho_Q(p^\nu)\ll \nu p^{2(t-\delta)\nu}$$
with $\delta:=\sum_{j=1}^t 1/\ell_j$.}
\smallskip
\nid The second estimate follows from \ref{ZS}. The first  formula  is included in the first estimate of \citeplus{R11}{lemma~3.4}, itself resting on an estimate of Korobov \citeplus{Ko92}{ch.\thinspace1, th.\thinspace5}.  The third bound is a consequence of the second estimate in \citeplus{R11}{lemma~3.4}, resting on Korobov's \citeplus{Ko92}{ch.\thinspace1, th.\thinspace6}, observing that the proof given in \citer{R11} remains valid as it stands when $\delta<1/2$. \qed
\medskip
We are now ready to complete the proof of \ref{majV2}.  Put $L:=\max \ell_j$. \ref{majV2} will follow from \ref{th}, applied with $k=1$, $Q_1=Q$,  $F=\Delta$.  Up to dividing through by $\gcd(c_1,\ldots,c_t)$, the $c_j$ may be assumed coprime. The polynomial $Q$ is then irreducible in $\z[X_1,\ldots, X_t, Y_1,\ldots, Y_t]$, so that $r=1$.  Indeed, we may write $Q=c_1X_1^{\ell_1}+Q_1$, where $Q_1\in \z[X_2,\ldots  X_t, Y_1,\ldots, Y_t]$ has no square factor in  $\Q[X_2,\ldots, X_t, Y_1,\ldots, Y_t]$. This enables applying Eisenstein's criterion associated to a non constant irreducible factor of $Q_1$ in $\z[X_2,\ldots, X_t, Y_1,\ldots, Y_t]$.
\par 
Taking the above into account, \ref{lemmeRobert} furnishes  $$
 \prod_{L< p\leqslant x}\Big(
1-{\varrho_Q^+(p)\over
p^{2t}}\Big)\asymp \log x\qquad (x\geqslant 2L). $$ Moreover,  the quantity
 $ E_Q(v)$ defined in \eqref{defER} satisfies
 $$ E_Q(v) \asymp \sum_{n\leqslant v} {\Delta(n)\over n}\qquad (v\geqslant 1). \eqdef{estEQ}$$
  This follows from Lemmas   \ref{slemma2} and \ref{slemmeRobert}. Indeed, since $k=r=1$, we have
  $$ E_Q(v)=\sum_{n\leqslant v} {\Delta(n)\rho_Q^{\#}(n)\over \K(n)^{2t}}\cdot\eqdef{EQD}$$
By \ref{lemmeRobert}, there exist two multiplicative functions  $\theta_1$ and $\theta_2$ such that 
$$ \eqalign{\rm{(i)}&\quad \theta_j(p^\nu)=1+O(1/p)\quad (j\in \{ 1,2\},\ \nu\geqslant 1),\cr
\rm{(ii)}&\quad \rho^{\#}(n)/\K(n)^{2t}=
\rho^{+}(n){\theta_1}(n)/n^{2t}=
\theta_2(n)/n\quad(n\geqslant 1).\cr}$$
A first application of \eqref{maj1sumH}
provides the upper bound in \eqref{estEQ}. To establish the corresponding lower bound, we note that, by condition (i), we have $\vartheta_2(p)\neq0$ for all but only a bounded number of primes $p$. Let $\vartheta_3(p^\nu)=\vartheta_2(p^\nu)+\1_{\{\vartheta_2(p^\nu)=0\}}(p^\nu)$. Applying \eqref{maj2sumH}
to replace  condition $n\leqslant v$ by $p|n\Rightarrow p\leqslant v$, we hence see that substituting $\vartheta_3$  to $\vartheta_2$ in \eqref{EQD}   does not alter the order of magnitude of the sum.  A second application of \eqref{maj1sumH} %\citeplus{H12}{lemma 2} 
with $H(n):=
\Delta(n)\theta_3(n)$ and $\theta=1/\theta_3$ then furnishes the lower bound in  \eqref{estEQ}.
\par 
We may now apply \ref{th} in order to derive \ref{majV2}.   
Put $$(\m,\n)=(m_1,\ldots, m_t, n_1,\ldots, n_t).$$ By symmetry, we may assume $n_0>m_0$. Let us then split the variation range of  $m_0+n_0$ into dyadic intervals $]N,2N]$ with $N=2^\nu\leqslant 2x^{1/\ell_0}
$. Up to neglecting a contribution 
$\ll x^{2\delta}$, %\citeplus{R11}{lemma~3.6} permits to 
 we may  assume  $N> x^{1/\ell_0-1/2\ell_0\ell_k}$. 
Indeed, if $N\leqslant  x^{1/\ell_0-1/2\ell_0\ell_k}$, then for all fixed $(m_1,\ldots,m_k)$,  $(n_1,\ldots,n_{k-1})$, the power $n_k^{\ell_k}$ lies in a fixed interval of length $\ll N^{\ell_0}\ll x^{1-1/2\ell_k}.$ The number of admissible $n_k$ is hence $\ll x^{1/\ell_k-1/2\ell_k^2}.$ However, for each fixed $n_k$, the number of admissible pairs $(m_0,n_0)$ is $\ll x^{\varepsilon}$ for any $\varepsilon>0$. This plainly implies that the contribution of   $N\leqslant  x^{1/\ell_0-1/2\ell_0\ell_k}$ is 
$ \ll x^{2\delta-1/2\ell_k^2+\varepsilon}\ll x^{2\delta}$ provided that  $\varepsilon\leqslant 1/2\ell_k^2.$

If $(m_0,\m,n_0,\n)$ is counted in $V_2^{\neq}(x)$, we have $Q(\m,\n)=c_0(n_0^{\ell_0}-m_0^{\ell_0})$,  hence the divisor $n_0^{\ell_0-1}+m_0n_0^{\ell_0-1}+\ldots +m_0^{\ell_0-1}$ of $Q(\m,\n)$ is of size $ N^{\ell_0-1}$. As a consequence, for each $2t$-tuple $(\m,\n)$, the number of admissible pairs $(m_0,n_0)$ does not exceed $\Delta( Q(\m,\n)).$
Now let us split the variation ranges of the $m_j^{\ell_j}$ and $n_h^{\ell_h}$ into intervals of size $N^{\ell_0}$, say  
$ m_j^{\ell_j}\in ]r_jN^{\ell_0},  (r_j+1)N^{\ell_0}] $ and $ n_h^{\ell_h}\in ]s_hN^{\ell_0},  (s_h+1)N^{\ell_0}] $ with $r_j,s_h\ll x/N^{\ell_0}$.
For each $r_j$, the integer $m_j$ is hence restricted to an interval of length $\ll N^{{\ell_0}/\ell_j}/(r_j+1)^{(1-1/\ell_j)}.$ 
\par 
For each pair  $(r_j,s_j)$, \ref{th} furnishes a contribution $$\ll
\prod_{ 1\leqslant j\leqslant t} {N^{2{\ell_0}/\ell_j } \over  (r_j+1)^{1-1/\ell_j}(1+s_j)^{1-1/\ell_j}}  {E_Q(x)\over \log x}\cdot\eqdef{majrs}$$ Here, the hypothesis $N> x^{1/\ell_0-1/2\ell_0\ell_k}$ has been used  in a crucial way   since  it implies that $m_j$ lies in a range of  size  $\gg x^{1/\ell_j-1/2\ell_k}\gg x^{ 1/2\ell_k} $. 
This enables us to check hypotheses \eqref{condth} of \ref{th}, which require that conditions \eqref{condth} be satisfied by
$$\eqalign{ &x_j=(r_j N^{\ell_0})^{1/\ell_j},\quad x_{j+k}=(s_j N^{\ell_0})^{1/\ell_j},\quad y_j=\{(r_j +1)N^{\ell_0}\}^{1/\ell_j}-(r_j N^{\ell_0})^{1/\ell_j},\cr
&y_{j+k}= ((s_j +1)N^{\ell_0})^{1/\ell_j}-(s_j N^{\ell_0})^{1/\ell_j},\cr}$$ for $1\leqslant j\leqslant k.$
Thus, the choice $\alpha=1/2\ell_k$,  $\beta=1/\ell_k-1/2\ell_k^2$, is admissible in \eqref{condth}.

Now observe that relations $\big|m_j^{\ell_j}-   n_j^{\ell_j}-(r_j-s_j)N^{\ell_0}\big|\leqslant N^{\ell_0}$ and
$$c_0|n_0^{\ell_0}-m_0^{\ell_0}|=\Big| \sum_{1\leqslant j\leqslant t}c_j (m_j^{\ell_j}-   n_j^{\ell_j})\Big|
\leqslant 2^{\ell_0}c_0N^{\ell_0}$$
imply
$
\Big| \sum_{1\leqslant j\leqslant t}c_j(r_j-  s_j)\Big|\leqslant  2^{\ell_0}c_0+\sum_{1\leqslant j\leqslant t}c_j.$ Therefore, when all $r_j$ and all $s_h$ are fixed except $s_1$,  the number of admissible values for  $s_1$
is  bounded.  Summing \eqref{majrs} over $r_j$ and  $s_h$, taking the condition  $\sum_{1\leqslant j\leqslant t} 1/\ell_j=\delta$ into account then yields a contribution 
$$\ll N^{\ell_0/\ell_1} x^{2\delta-1/\ell_1} {E_Q(x)\over \log x}\cdot$$
It remains to sum over $\nu$ such that $N^{\ell_0}=2^{\ell_0\nu}\leqslant 2^{\ell_0}x $ to get, in view of \eqref{estEQ}, 
$$V_2^{\neq }(x;\cc,\,\el)\ll x^{2\delta}{E_Q(x)\over \log x}\ll x^{2\delta}\gS(x).$$
Note that the manipulation resting on the condition $m_0+n_0\in ]N,2N]$  saved a factor  $\log_2 x$. \par 
This completes the proof of \ref{majV2}.\par
\medskip 

\paradeuxb{Proof of \ref{thRobert}.}
Since $V_1(x)\asymp x$ (see \citeplus{R11}{lemma 3.1} for further precisions) and taking \eqref{inegCSV} into account  the required estimate \eqref{estV0} immediately results from a suitable bound for $V_2(x,\cc,\el)$.  In view of \eqref{estV2neq}, we write $\cc=(c_1,\ldots,c_t),\,\el=(\ell_1,\ldots,\ell_t)$.\par 
 For $1\leqslant k\leqslant t$, put $\delta_k=\sum_{t-k\leqslant j\leqslant t} 1/\ell_j $,  $\cc_k=(c_{t-k},\ldots, c_t)$, $\el_k=(\ell_{t-k},\ldots, \ell_t)$. The parameter $s$ being defined in the statement of \ref{thRobert} and $\el$ satisfying conditions (i), (ii), (iii) of this statement, we shall show by induction  on $k\geqslant \max(1,s-1)$ that 
$$V_2(x;\cc_k,\el_k)\ll x^{\delta_k}\gS(x) .  \eqdef{estV2finale}$$
  Note that hypothesis  (iii) may be rewritten as 
$  \delta_k\leqslant 1/\ell_{t-k-1}$ for $s -1 \leqslant k\leqslant t -1. $

 For  $k=s=1$, $t\geqslant 2$, the stated bound \eqref{estV2finale}     follows from \citeplus{R11}{lemma 3.11}. The case  $k=1$, $s=2$, $t\geqslant 2$, follows from \eqref{estV2check}. For  $k=2$, $t\geqslant 3$, $s=3$,  we only need to bound the number of solutions $(\x,\y)$ such that $x_j, y_j\leqslant x^{1/\ell}$ and 
$$x_1^\ell+x_2^\ell+x_3^\ell=y_1^\ell+y_2^\ell+y_3^\ell $$ with $\ell=\ell_t$.
The required estimate then follows from \ref{majV2=}.
The induction is thus initialized. In these cases, the factor $\gS(x) $ is actually superfluous.
\par 
Let $k\geqslant s-1$ and assume \eqref{estV2finale} holds for $k$. We want to prove that it persists for $k+1$.
By Propositions \ref{smajV2=} and \ref{smajV2}, we have
$$V_2(x ;\cc_{k+1},\,\el_{k+1})\ll x^{ 2\delta_{k}}\gS(x)  +x^{1/\ell_{t-k-1}} V_2(x ;\cc_{k},\el_{k})
\ll (x^{ 2\delta_{k} }  +x^{1/\ell_{t-k-1}+\delta_k})\gS(x) .$$
This provides the desired bound  since $\delta_k\leqslant 1/\ell_{t-k-1}$, and $1/\ell_{t-k-1}+\delta_k= \delta_{k+1}$.
\par 
The induction step is thus established, and so we get the result for $k=t+1$. This completes the proof of \ref{thRobert} since $\delta_{t+1}=1$.

\medskip\goodbreak
\goodbreak

\medskip
\noi{\bf Acknowledgement.} The authors take pleasure in addressing warm thanks to Tim Browning, Carlo Pagano, Per Salberger,  and the referee for  pertinent  remarks and suggestions. 

\goodbreak
 
\bigskip
 \centerline{\twelvebf References}
\bigskip
 {\eightpoint\leftskip9mm\rightskip5mm
 
 \bibtem{BB95} A. Balog \& J. BrŸdern, Sums of three cubes in three linked three-progressions, {\it J. reine und angew. Math. \bf 466} (1995), 45Ð85.\par 
\bibtem{dlBB06} R. de la Bretche \& T.D.
Browning,  {Sums of arithmetic functions over
values of binary forms},    {\it Acta
Arith.} {\bf 125} n${}^{\circ}$  3 (2006), 291-304.
\par
\bibtem{dlBT12} R. de la Bretche \&  G. Tenenbaum, Moyennes de fonctions arithmŽtiques de formes binaires,
   {\it Mathematika}   {\bf 58} (2012), 290-304.
   \par

\bibtem{dlBT24} R. de la Bretche \&  G.
Tenenbaum,
Note on the mean value of the Erd\H os-Hooley Delta function, {\it Acta Arith. \bf 219.4} (2015), 379-394. 
\par

\bibtem{BHB04} T.D. Browning \& D.R. Heath-Brown,	Equal sums of three powers, {\it Invent. Math.}Ê{\bf 157} (2004), 553--573.
\par

\bibtem{CKPS24} S. Chan, P. Koymans, C. Pagano \& E. Sofos,
6-torsion and integral points on quartic surfaces, Ann. Sc. Norm. Super. Pisa Cl. Sci., to appear.
\par
 \bibtem{FGK23} K. Ford,  B. Green, D. Koukoulopoulos, Equal sums in random sets and the concentration of divisors, 
{\it Inventiones Math.} {\bf 232}, (2023), 1027--1160. \par 
\bibtem{Er73} P. Erd\H os, Problem 218, Solution by the proposer, {\it Can. Math. Bull. \bf 16} (1973),  621-622.\par 
\bibtem{FKT23} K. Ford, D. Koukoulopoulos, T. Tao, A lower bound on the mean value of the Erd\H os--Hooley Delta function,
{\it Proc. Lond. Math. Soc.} (3) {\bf 129} (2024), no. 1, Paper No. e12618, 18 pp.
% preprint, https://arxiv.org/abs/2308.11987.
\par
\bibtem{HR74}
H. Halberstam, H.-E. Richert, {\it Sieve Methods}, London Math. Soc. Monogr., no. 4,
Academic Press, London-New York, 1974, xiv+364 pp.
\par
\bibtem{HT82} R.R. Hall \& G. Tenenbaum, On the average and normal orders of Hooley's $\Delta$-function,
{\it J. London Math. Soc.} (2) {\bf 25} (1982), 392-406. 
\par 
\bibtem{HT88} R.R. Hall and G. Tenenbaum, {\it Divisors}, Cambridge tracts in
mathematics 90, Cambridge University Press (1988, paperback ed. 2008).

\bibtem{H12} K. Henriot, NairÐTenenbaum bounds uniform with respect to the discriminant, {\it 
    Math. Proc. Camb. Phil. Soc. \bf 152}, no. 3 (2012),  405--424.
 \par    
    
 \bibtem{H14} K. Henriot,    NairÐTenenbaum uniform with respect to the discriminant--Erratum, {\it  Math. Proc. Camb. Phil. Soc. \bf 157} (2014), 375--377.
\par     
\bibtem{Ho79} C. Hooley,  A new technique and its
applications to the theory of numbers, {\it Proc.
London Math.
Soc.} (3) {\bf 38} (1979), 115--151.
\bibtem{Ko92} N.M. Korobov, {\it Exponential sums and their applications}, Math. Appl. (Soviet Ser.), 80,
Kluwer Academic Publishers Group, Dordrecht, 1992, xvi+209 pp.\par 
 \bibtem{KT23} D. Koukoulopoulos \& T. Tao, An upper bound on the mean value of the Erd\H os--Hooley Delta function, {\it Proc. London Math. Soc. (3) \bf 127}, no.6 (2023),  1865Ð1885.

\bibtem{Na92} M. Nair, Multiplicative functions
of polynomial values in short intervals,
{\it Acta Arith. \bf 62} (1992), 257--269.\par
\bibtem{NT98} M. Nair \& G. Tenenbaum, Short sums
of certain arithmetic functions, {\it Acta Math.}
{\bf 180}, (1998), 119-144.

\bibtem{PSW16} L.B. Pierce, D. Schindler, M. Wood, Representations of integers by systems of three quadratic
forms, {\it Proc. Lond. Math. Soc.} (3) {\bf 113} (2016), 289--344.\par 

\bibtem{R11} O. Robert,
Sur le nombre des entiers reprŽsentables comme somme de trois puissances
{\it Acta Arith.} {\bf 149} (2011), 1--21.\par

\bibtem{Sa05}  P. Salberger, Counting rational points on hypersurfaces of low dimension, {\it Ann. Sci. ƒc. Norm. Sup.} {\bf 38} (2005), 93--115.
\par 

\bibtem{Sa23}  P. Salberger, Counting rational points on projective varieties,  
{\it Proc. London Math. Soc.}
{\bf 126},  4 (2023), 1092--1133.
\par 

\bibtem{Sa24}  P. Salberger, Private communication, April 2024.
\par

\bibtem{S91}  C.L. Stewart, On the number of
solutions of polynomial congruences and Thue
equations, {\it J. Amer. Math. Soc.}  {\bf 4}
(1991), 793--835.

\bibtem{Sh80} P. Shiu, A Brun--Titchmarsh theorem
for multiplicative functions, {\it J. reine
angew. Math. \bf313} (1980), 161--170.

\bibtem{S80} 
J.T. Schwartz,   Fast Probabilistic Algorithms for Verification of Polynomial Identities, {\it J. Assoc. Comput. Math.} {\bf27}, no 4, (1980), 701--717.

\bibtem{Te85} G. Tenenbaum, Sur la concentration
moyenne des diviseurs,  {\it Comment. Math.
Helvetici \bf 60}
(1985), 411-428.\par
\bibtem{Te86} G. Tenenbaum,
Fonctions $\Delta$ de Hooley et applications, SŽminaire de thŽorie des nombres, Paris 1984-85, Prog. Math. 63 (1986), 225-239.
\par
\bibtem{GT15} G. Tenenbaum, {\it Introduction to analytic and probabilistic number theory}, Graduate Studies in Mathematics 163, Amer. Math. Soc. 2015; see also {\it Introduction ˆ la thŽorie analytique et probabiliste des
nombres}, 5th ed., Dunod, Sciences Sup, 2022, 547 pp.
\par 
\bibtem{W18} V.Y. Wang, 
Sums of cubes and the Ratios Conjectures, preprint, https://arxiv.org/abs/2108. 03398.
 \par
\bibtem{Z79} 
R. Zippel,  Probabilistic algorithms for sparse polynomials, International Symposium on Symbolic and Algebraic Manipulation (EUROSAM) 1979, {\it Symbolic and Algebraic Computation}, 1979, p. 216--226.\par }

{\leftskip2mm\rightskip-2cm\sevenrm
\gutter=4cm \doublecolumns
 \obeylines \baselineskip=7pt
RŽgis de la Bretche
UniversitŽ Paris CitŽ, Sorbonne UniversitŽ, CNRS,
Institut Universitaire de France,
Institut de Math. de Jussieu-Paris Rive Gauche
 F-75013 Paris  
France
\smallskip
{\seventt regis.delabreteche@imj-prg.fr}
\phantom a\goodbreak
 GŽrald Tenenbaum
Institut \'Elie Cartan
Universit\'e de Lorraine
 BP 70239\par
54506 Vand\oe uvre-ls-Nancy Cedex
 France
\smallskip
{\seventt gerald.tenenbaum@univ-lorraine.fr}
\singlecolumn
\par}

\vfill\eject

\end